\documentclass{article}
\usepackage{amssymb}

\input{tcilatex}
\begin{document}

\begin{center}
\bigskip \textbf{Multipoint Cauchy problem for nonlinear wave equations in
vector-valued spaces}

\textbf{Veli Shakhmurov}

Department of Mechanical Engineering, Okan University, Akfirat, Tuzla 34959
Istanbul, Turkey,

E-mail: veli.sahmurov@okan.edu.tr

A\textbf{bstract}
\end{center}

In this paper, regularity properties, Strichartz type estimates to solutions
of multipo\i nt Cauchy problem for linear and nonlinear abstract wave
equations in vector-valued function spaces are obtained. The equation
includes a linear operator $A$ defined in a Hilbert space $H$, in which by
choosing $H$ and $A$ we can obtain numerous classis of nonlocal initial
value problems for wave equations which occur in a wide variety of physical
systems.


\begin{center}
\ \ \textbf{AMS 2010: 35Q41, 35K15, 47B25, 47Dxx, 46E40 }

\textbf{1. Introduction, definitions}
\end{center}

\bigskip\ Consider the multipoint Cauchy problem for nonlinear abstract wave
equations (NLAWE)

\begin{equation}
\partial _{t}^{2}u-\Delta u+Au=F\left( u\right) ,\text{ }x\in R^{n},\text{ }%
t\in \left[ 0,T\right] ,  \tag{1.1}
\end{equation}%
\begin{equation}
u\left( 0,x\right) =\varphi \left( x\right) +\dsum\limits_{k=1}^{m}\alpha
_{k}u\left( \lambda _{k},x\right) ,\text{ for a.e. }x\in R^{n},  \tag{1.2}
\end{equation}%
\begin{equation}
u_{t}\left( 0,x\right) =\psi \left( x\right) +\dsum\limits_{k=1}^{m}\beta
_{k}u_{t}\left( \lambda _{k},x\right) ,\text{ for a.e. }x\in R^{n}, 
\tag{1.3}
\end{equation}%
where $m$ is a positive integer, $\alpha _{k}$, $\beta _{k}$ are complex
numbers, $\lambda _{k}\in \left( 0,\right. \left. T\right] ,$ $A$ is a
linear and $F$ is a nonlinear operators in a Hilbert space $E$, $\Delta $
denotes the Laplace operator in $R^{n}$, $u=$ $u(t,x)$ is the $E$-valued
unknown function and $\varphi $, $\psi $ are data functions.

If we put $F\left( u\right) =$ $\lambda \left\vert u\right\vert ^{p}u$ in $%
\left( 1.1\right) $ we get the multipoint initial value problem for the
following NLAWE%
\begin{equation}
\partial _{t}^{2}u-\Delta u+Au=\lambda \left\vert u\right\vert ^{p}u,\text{ }%
x\in R^{n},\text{ }t\in \left[ 0,T\right] ,  \tag{1.4}
\end{equation}%
\[
u\left( 0,x\right) =\varphi \left( x\right) +\dsum\limits_{k=1}^{m}\alpha
_{k}u\left( \lambda _{k},x\right) ,\text{ for a.e. }x\in R^{n}, 
\]%
\[
u_{t}\left( 0,x\right) =\psi \left( x\right) +\dsum\limits_{k=1}^{m}\beta
_{k}u_{t}\left( \lambda _{k},x\right) ,\text{ for a.e. }x\in R^{n}, 
\]%
where $p\in \left( 1,\infty \right) $, $\lambda $ is a real number.

Let $\mathbb{N}$, $\mathbb{R}$ and $\mathbb{C}$\ denote the sets of all
natural, real and complex numbers, respectively. For $E=\mathbb{C}$, $\alpha
_{k}=\beta _{k}=0$ and $A=0$ the problem $\left( 1.4\right) $\ become the
classical Cauchy problem for nonlinear wave equation (NWE)%
\begin{equation}
\partial _{t}^{2}u-\Delta u=\lambda \left\vert u\right\vert ^{p-1}u,\text{ }%
x\in R^{n},\text{ }t\in \left[ 0,T\right] ,  \tag{1.5}
\end{equation}%
\[
u\left( 0,x\right) =\varphi \left( x\right) \text{, }u_{t}\left( 0,x\right)
=\psi \left( x\right) \text{ for a.e. }x\in R^{n}. 
\]

The existence of solutions and regularity properties of Cauchy problem for
NWE studied e.g in $\left[ 4\right] $, $\left[ 6\right] $, $\left[ 9\right] $%
, $\left[ 14\right] $, $\left[ 17\right] ,$ $\left[ 19-22\right] ,$ $\left[
24\right] ,$ $\left[ 27\right] ,$ $\left[ 30-31\right] ,$ $\left[ 35\text{, }%
37\right] ,$ $\left[ 39-42\right] $ and the referances therein.\ In
contrast, to the mentioned above results we will study the regularity
properties of the abstract Cauchy problem $\left( 1.1\right) $. Abstract
differential equations studied e.g. in $\left[ 1-3\right] $, $\left[ 7-8%
\right] $, $\left[ 11\right] $, $\left[ 13\right] $, $\left[ 15-16\right] $, 
$\left[ 18\right] $, $\left[ 25\right] $, $\left[ 28-29\right] $, $\left[
32-34\right] ,$ $\left[ 41\right] $ and $\left[ 43\right] .$ Since the
Hilbert space $H$ is arbitrary and $A$ is a possible linear operator, by
choosing $H$ and $A$ we can obtain numerous classes of wave equations and
its systems which occur in a wide variety of physical systems.

Our main goal is to obtain the exsistence, uniquness and Strichartz type
estimates, i.e. estimates in the form of space time integrability to
solution of $\left( 1.1\right) -\left( 1.3\right) .$ Strichartz type
estimates to solutions of Cauchy problem for wave equations studied e.g in $%
\left[ 10\right] $, $\left[ 14\right] $, $\left[ 20\right] $, $\left[ 22%
\right] ,$ $\left[ 24\right] ,$ $\left[ 35\right] ,$ $\left[ 37\right] .$ If
we\ choose $H$ a concrete space, for example $H=L^{2}\left( \Omega \right) $%
, $A=L,$ where $\Omega $ is a domin in $R^{d}$ with sufficiently smooth
boundary, in variables $y=\left( y_{1},y_{2},...y_{d}\right) $ and $L$ is an
elliptic operator in $L^{2}\left( \Omega \right) $ in $\left( 1.2\right) ,$
then we obtain exsistence, uniquness and the regularity properties of the
mixed problem for linear wave equation

\begin{equation}
\partial _{t}^{2}u-\Delta u+Lu=F\left( t,x,y\right) ,\text{ }t\in \left[ 0,T%
\right] \text{, }x\in R^{n},\text{ }y\in \Omega ,  \tag{1.6}
\end{equation}%
and for the NLS equation%
\[
\partial _{t}^{2}u-\Delta u+Lu=F\left( u\right) ,\text{ }t\in \left[ 0,T%
\right] \text{, }x\in R^{n},\text{ }y\in \Omega ,
\]%
where $u=u\left( t,x,y\right) .$

\ Moreover, let we choose $E=L^{2}\left( 0,1\right) $ and $A$ to be
differential operator with generalized Wentzell-Robin boundary condition
defined by 
\[
D\left( A\right) =\left\{ u\in W^{2,2}\left( 0,1\right) ,\text{ }%
B_{j}u=Au\left( j\right) =0,\text{ }j=0,1\right\} ,\text{ } 
\]%
\[
\text{ }Au=au^{\left( 2\right) }+bu^{\left( 1\right) } 
\]%
where $a=a\left( y\right) $ and $b=b\left( y\right) $ are complex-valued
functions. Then, from the main our theorem we get the exsistence, uniqueness
and regularity properties of multipoint Wentzell-Robin type mixed problem
for the following wave equation 
\begin{equation}
\partial _{t}^{2}u-\Delta u+a\frac{\partial ^{2}u}{\partial y^{2}}+b\frac{%
\partial u}{\partial y}=F\left( t,x\right) ,\text{ }  \tag{1.7}
\end{equation}%
\ \ \ 

\begin{equation}
B_{j}u=0\text{, }j=0,1.  \tag{1.8}
\end{equation}

\begin{equation}
u\left( 0,x,y\right) =\varphi \left( x,y\right)
+\dsum\limits_{k=1}^{m}\alpha _{k}u\left( \lambda _{k},x,y\right) ,\text{
for a.e. }x\in R^{n},\text{ }y\in \left( 0,1\right) ,  \tag{1.9}
\end{equation}%
\[
u_{t}\left( 0,x,y\right) =\psi \left( x,y\right)
+\dsum\limits_{k=1}^{m}\beta _{k}u_{t}\left( \lambda _{k},x,y\right) ,\text{
for a.e. }x\in R^{n}\text{, }y\in \left( 0,1\right) 
\]%
and the same mixed problem for the following NWE%
\begin{equation}
\partial _{t}^{2}u-\Delta u+a\frac{\partial ^{2}u}{\partial y^{2}}+b\frac{%
\partial u}{\partial y}=F\left( u\right) ,\text{ }  \tag{1.10}
\end{equation}%
where 
\[
u=u\left( t,x,y\right) \text{, }t\in \left[ 0,T\right] \text{, }x\in R^{n},%
\text{ }y\in \left( 0,1\right) . 
\]%
Note that, the regularity properties of Wentzell-Robin type BVP for elliptic
equations were studied e.g. in $\left[ \text{12, 23 }\right] $ and the
references therein. Moreover, if put $E=l_{2}$ and choose $A$ as a infinite
matrix $\left[ a_{mj}\right] $, $m,j=1,2,...,\infty ,$ then from our results
we obtain the exsistence, uniquness and regularity properties of multipoint
Cauchy problem for infinity many system of linear wave equations 
\begin{equation}
\partial _{t}^{2}u_{m}-\Delta
u_{m}+\sum\limits_{j=1}^{N}a_{mj}u_{j}=F_{j}\left( t,x\right) ,\text{ }t\in %
\left[ 0,T\right] \text{, }x\in R^{n},  \tag{1.11}
\end{equation}%
\[
u_{m}\left( 0,x\right) =\varphi _{m}\left( x\right)
+\dsum\limits_{k=1}^{m}\alpha _{k}u_{m}\left( \lambda _{k},x\right) ,\text{
for a.e. }x\in R^{n}, 
\]%
\[
\partial _{t}u_{m}\left( 0,x\right) =\psi _{m}\left( x\right)
+\dsum\limits_{k=1}^{m}\beta _{k}\partial _{t}u_{m}\left( \lambda
_{k},x\right) ,\text{ for a.e. }x\in R^{n} 
\]%
and the same problem for infinity many system of NWE equation 
\begin{equation}
\partial _{t}^{2}u_{m}-\Delta
u_{m}+\sum\limits_{j=1}^{N}a_{mj}u_{j}=F_{m}\left(
u_{1},u_{2},...u_{N}\right) ,\text{ }t\in \left[ 0,T\right] \text{, }x\in
R^{n},  \tag{1.12}
\end{equation}%
where $a_{mj}$ are complex numbers, $u_{j}=u_{j}\left( t,x\right) .$

\begin{center}
\textbf{2. Definitions and} \textbf{Background}
\end{center}

Let $E$ be a Banach space. $L^{p}\left( \Omega ;E\right) $ denotes the space
of strongly measurable $E$-valued functions that are defined on the
measurable subset $\Omega \subset R^{n}$ with the norm

\[
\left\Vert f\right\Vert _{p}=\left\Vert f\right\Vert _{L^{p}\left( \Omega
;E\right) }=\left( \int\limits_{\Omega }\left\Vert f\left( x\right)
\right\Vert _{E}^{p}dx\right) ^{\frac{1}{p}},\text{ }1\leq p<\infty \ . 
\]

Let $H$ be a Hilbert space. For $p=2$ and $E=H$ the space $L^{p}\left(
\Omega ;E\right) $ become the $H$-valued functions space $L^{2}\left( \Omega
;H\right) $ with inner product:%
\[
\left( f,g\right) _{L^{2}\left( \Omega ;H\right) }=\int\limits_{\Omega
}\left( f\left( x\right) ,g\left( x\right) \right) _{H}dx\text{, for any }%
f,g\in L^{2}\left( \Omega ;H\right) . 
\]

Let $L_{t}^{q}L_{x}^{r}\left( E\right) =L_{t}^{q}L_{x}^{r}\left( \left(
a,b\right) \times \Omega ;E\right) $ denotes the space of strongly
measurable $E$-valued functions that are defined on the measurable set $%
\left( a,b\right) \times \Omega $ with the norm 
\[
\left\Vert f\right\Vert _{L_{t}^{q}L_{x}^{r}\left( \left( a,b\right) \times
\Omega ;E\right) }=\left( \dint\limits_{a}^{b}\left[ \int\limits_{\Omega
}\left\Vert f\left( t,x\right) \right\Vert _{E}^{r}dx\right] ^{\frac{q}{r}%
}dt\right) ^{\frac{1}{q}},\text{ }1\leq q,r<\infty \ . 
\]

Let $C\left( \Omega ;E\right) $ denote the space of $E-$valued, bounded
strongly continious functions on $\Omega $ with norm 
\[
\left\Vert u\right\Vert _{C\left( \Omega ;E\right) }=\sup\limits_{x\in
\Omega }\left\Vert u\left( x\right) \right\Vert _{E}. 
\]

$C^{m}\left( \Omega ;E\right) $\ will denote the spaces of $E$-valued
bounded strongly continuous and $m$-times continuously differentiable
functions on $\Omega $ with norm 
\[
\left\Vert u\right\Vert _{C^{m}\left( \Omega ;E\right) }=\max\limits_{0\leq
\left\vert \alpha \right\vert \leq m}\sup\limits_{x\in \Omega }\left\Vert
D^{\alpha }u\left( x\right) \right\Vert _{E}. 
\]

Let $E_{1}$ and $E_{2}$ be two Banach spaces. $B\left( E_{1},E_{2}\right) $
will denote the space of all bounded linear operators from $E_{1}$ to $%
E_{2}. $ For $E_{1}=E_{2}=E$ it will be denoted by $B\left( E\right) .$

A closed densely defined linear operator\ $A$ is said to be absulute
positive in a Banach\ space $E$ (see $\left[ 11\right] $, \S\ 11.2) if $%
D\left( A\right) $ is dense on $E,$ the resolvent $\left( A-\lambda
^{2}I\right) ^{-1}$ exists for $\func{Re}\lambda >\omega $ and 
\[
\left\Vert \left( A-\lambda ^{2}I\right) ^{-1}\right\Vert _{B\left( E\right)
}\leq M_{0}\left\vert \func{Re}\lambda -\omega \right\vert ^{-1}\text{. } 
\]

It is known $\left[ \text{38, \S 1.15.1}\right] $ that there exist
fractional powers\ $A^{\theta }$ of a absolute positive operator $A.$ Let $%
E\left( A^{\theta }\right) $ denote the space $D\left( A^{\theta }\right) $
with the graphical norm 
\[
\left\Vert u\right\Vert _{E\left( A^{\theta }\right) }=\left( \left\Vert
u\right\Vert ^{p}+\left\Vert A^{\theta }u\right\Vert ^{p}\right) ^{\frac{1}{p%
}},1\leq p<\infty ,\text{ }0<\theta <\infty . 
\]

For case of Hilber space $H$ and $p=2$, $E\left( A^{\theta }\right) $ will
be denoted by $H\left( A^{\theta }\right) .$

\textbf{Remark 1.1. }It is known\ that if the operator $A$ is absolute
positive in a Banach\ space $E$ and $0\leq \alpha <1$ then it is an
infinitesimal generator of group of bounded linear operator $U_{A}\left(
t\right) $ satisfying 
\[
\left\Vert U_{A}\left( t\right) \right\Vert _{B\left( E\right) }\leq
Me^{\omega \left\vert t\right\vert },\text{ }t\in \left( -\infty ,\infty
\right) , 
\]%
\begin{equation}
\left\Vert A^{\alpha }U_{A}\left( t\right) \right\Vert _{B\left( E\right)
}\leq M\left\vert t\right\vert ^{-\alpha },\text{ }t\in \left( -\infty
,\infty \right)  \tag{2.1}
\end{equation}%
(see e.g. $\left[ \text{29, \S\ 1.6}\right] $, Theorem 6.3).

Let $E$ be a Banach space. $S=S(R^{n};E)$ denotes $E$-valed Schwartz class,
i.e. the space of all $E$ -valued rapidly decreasing smooth functions on $%
R^{n}$ equipped with its usual topology generated by seminorms. $S(R^{n};%
\mathbb{C})$ denoted by $S$.

Let $S^{\prime }(R^{n};E)$ denote the space of all continuous linear
operators, $L:S\rightarrow E$, equipped with the bounded convergence
topology. Recall $S(R^{n};E)$ is norm dense in $L^{p}(R^{n};E)$ when $%
1<p<\infty .$

Let $F$ denotes the Fourier trasformation, $\hat{u}=Fu$ and 
\[
s\in \mathbb{R},\text{ }\xi =\left( \xi _{1},\xi _{2},...,\xi _{n}\right)
\in R^{n},\text{ }\left\vert \xi \right\vert ^{2}=\dsum\limits_{k=1}^{n}\xi
_{k}^{2}, 
\]%
\[
\langle \xi \rangle :=\left( 1+\left\vert \xi \right\vert ^{2}\right) ^{%
\frac{1}{2}}. 
\]%
Consider the $E-$valued Sobolev space $W^{s,p}\left( R^{n};E\right) $ and
homogeneous Sobolev spaces $\mathring{W}^{s,p}\left( R^{n};E\right) $\
defined by respectively, 
\[
W^{s,p}\left( R^{n};E\right) =\left\{ u:u\in S^{\prime }(R^{n};E),\right. 
\text{ } 
\]%
\[
\left\Vert u\right\Vert _{W^{s,p}\left( R^{n};E\right) }=\left. \left\Vert
F^{-1}\left( 1+\left\vert \xi \right\vert ^{2}\right) ^{\frac{s}{2}}\hat{u}%
\right\Vert _{L^{p}\left( R^{n};H\right) }<\infty \right\} , 
\]

\[
\mathring{W}^{s,p}\left( R^{n};E\right) =\left\{ u:u\in S^{\prime
}(R^{n};E),\left\Vert u\right\Vert _{\mathring{W}^{s,p}\left( R^{n};E\right)
}=\left\Vert F^{-1}\left\vert \xi \right\vert ^{s}\hat{u}\right\Vert
_{L^{p}\left( R^{n};H\right) }<\infty \right\} .\text{ } 
\]%
For $\Omega =R^{n}\times G,$ $\mathbf{p=}\left( p_{1},\text{ }p_{2}\right) ,$
$s\in \mathbb{R}$ and $l\in \mathbb{N}$ we define the $E$-valud anisotropic
Sobolev space $W^{s,l,p}\left( \Omega ;E\right) $ by 
\[
W^{s,l,\mathbf{p}}\left( \Omega ;E\right) :=\left\{ u\in S^{\prime }(\Omega
;E),\text{ }\left\Vert u\right\Vert _{W^{s,l,\mathbf{p}}\left( \Omega
\right) }=\left\Vert u\right\Vert _{W^{s,\mathbf{p}}\left( \Omega \right)
}+\left\Vert u\right\Vert _{W^{l,\mathbf{p}}\left( \Omega \right) }\right\}
, 
\]%
where 
\[
\left\Vert u\right\Vert _{W^{s,\mathbf{p}}\left( \Omega ;E\right)
}=\left\Vert F^{-1}\left( 1+\left\vert \xi \right\vert ^{2}\right) ^{\frac{s%
}{2}}\hat{u}\right\Vert _{L^{\mathbf{p}}\left( \Omega ;E\right) }<\infty , 
\]%
\[
\left\Vert u\right\Vert _{W^{l,\mathbf{p}}\left( \Omega ;E\right)
}=\left\Vert u\right\Vert _{L^{\mathbf{p}}\left( \Omega ;E\right)
}+\dsum\limits_{\left\vert \beta \right\vert =l}\left\Vert D_{y}^{\beta
}u\right\Vert _{L^{\mathbf{p}}\left( \Omega ;E\right) }. 
\]%
The similar way, we define homogeneous anisotropic Sobolev spaces $\mathring{%
W}^{s,l,\mathbf{p}}\left( \Omega ;E\right) $ as:%
\[
\mathring{W}^{s,l,\mathbf{p}}\left( \Omega ;E\right) :=\left\{ u\in
S^{\prime }(\Omega ;E),\text{ }\left\Vert u\right\Vert _{W^{s,l,\mathbf{p}%
}\left( \Omega ;E\right) }=\left\Vert u\right\Vert _{W^{s,\mathbf{p}}\left(
\Omega ;E\right) }+\left\Vert u\right\Vert _{\mathring{W}^{l,\mathbf{p}%
}\left( \Omega ;E\right) }\right\} , 
\]%
where 
\[
\left\Vert u\right\Vert _{\mathring{W}^{s,\mathbf{p}}\left( \Omega ;E\right)
}=\left\Vert F^{-1}\left\vert \xi \right\vert ^{s}\hat{u}\right\Vert _{L^{%
\mathbf{p}}\left( \Omega ;E\right) }<\infty . 
\]

Let $A$ be a linear operator in a Banach space $E.$ Consider Sobolev-Lions
type homogeneous and inhomogeneous abstract spaces, respectively 
\[
\mathring{W}^{s,p}\left( R^{n};E\left( A\right) ,E\right) =\mathring{W}%
^{s,p}\left( R^{n};E\right) \cap L^{p}\left( R^{n};E\left( A\right) \right) ,%
\text{ } 
\]%
\[
\left\Vert u\right\Vert _{\mathring{W}^{s,p}\left( R^{n};E\left( A\right)
,E\right) }=\left\Vert u\right\Vert _{\mathring{W}^{s,p}\left(
R^{n};E\right) }+\left\Vert u\right\Vert _{L^{p}\left( R^{n};E\left(
A\right) \right) }<\infty , 
\]

\[
W^{s,p}\left( R^{n};E\left( A\right) ,E\right) =W^{s,p}\left( R^{n};E\right)
\cap L^{p}\left( R^{n};E\left( A\right) \right) ,\text{ } 
\]%
\[
\left\Vert u\right\Vert _{W^{s,p}\left( R^{n};E\left( A\right) ,E\right)
}=\left\Vert u\right\Vert _{W^{s,p}\left( R^{n};E\right) }+\left\Vert
u\right\Vert _{L^{p}\left( R^{n};E\left( A\right) \right) }<\infty . 
\]

Sometimes we use one and the same symbol $C$ without distinction in order to
denote positive constants which may differ from each other even in a single
context. When we want to specify the dependence of such a constant on a
parameter, say $\alpha $, we write $C_{\alpha }$.

\textbf{Definition 1.1.} (Solution). A function $u$ : $\left[ 0,T\right]
\times R^{n}\rightarrow H\left( A\right) $ is called a (strong) solution to
problem $(1.1)-\left( 1.3\right) $ if it lies in the class 
\[
C_{t}^{0}\left( \left[ 0,T\right] ;\mathring{W}_{x}^{2,\gamma }\left(
R^{n};H\left( A\right) \right) \right) \cap C_{t}^{1}\left( \left[ 0,T\right]
;\mathring{W}_{x}^{2,\gamma -1}\left( R^{n};H\left( A\right) \right) \right) 
\]%
for $\gamma \geq n\left( \frac{1}{2}-\frac{1}{k}\right) $, $k>1$\ and obeys
the formula 
\[
u\left( t,x\right) =U_{\Delta +A}\left( t\right) \left[ \varphi \left(
x\right) +\dsum\limits_{k=1}^{m}\alpha _{k}u\left( \lambda _{k},x\right) %
\right] +
\]%
\begin{equation}
\tilde{U}_{\Delta +A}\left( t\right) \left[ \psi \left( x\right)
+\dsum\limits_{k=1}^{m}\beta _{k}u\left( \lambda _{k},x\right) \right]
+\dint\limits_{0}^{t}\tilde{U}_{\Delta +A}\left( t-s\right) F\left( u\left(
s\right) \right) ds\text{ }  \tag{2.2}
\end{equation}%
for all $t\in \left( 0,T\right) ,$ where $U_{A+\Delta }\left( t\right) $ is
a cosine, $\tilde{U}_{A+\Delta }\left( t\right) $ is a sine
operator-functions (see e.g. $\left[ 11\right] $) with generatr of $A+\Delta 
$, i.e.%
\begin{equation}
U_{A+\Delta }\left( t\right) =\frac{1}{2}\left( e^{t\left( A-\Delta \right)
^{\frac{1}{2}}}+e^{-t\left( A-\Lambda \right) ^{\frac{1}{2}}}\right) \text{, 
}  \tag{2.3}
\end{equation}%
\[
\tilde{U}_{A+\Delta }\left( t\right) =\frac{1}{2}A^{-\frac{1}{2}}\left(
e^{t\left( A-\Delta \right) ^{\frac{1}{2}}}-e^{-t\left( A-\Lambda \right) ^{%
\frac{1}{2}}}\right) .
\]

We say that $u$ is a global solution if $T=\infty $.

We write $a\lesssim b$ to indicate that $a\leq Cb$ for some constant $C$,
which is permitted to depend on some parameters.\bigskip

\begin{center}
\textbf{3. The exsistence of solution to multipoint Cauchy problem for
linear wave equation}
\end{center}

Consider the abstract Schr\"{o}dinger\ equation%
\begin{equation}
\partial _{t}^{2}u-\Delta u+Au=0,\text{ }t\in \left[ 0,T\right] ,\text{ }%
x\in R^{n},  \tag{3.1}
\end{equation}%
where $A$ is a linear operator in a Hilber space $H.$

It can be shown that the fundamental solutions of the free abstract Schr\"{o}%
dinger equation $\left( 3.1\right) $\ can be exspressed as 
\begin{equation}
U_{A+\Delta }\left( t\right) \left( x,y\right) =C\left( t,A\right) U_{\Delta
}\left( t\right) \left( x,y\right) \text{,}  \tag{3.2}
\end{equation}%
\[
\text{ }\tilde{U}_{A+\Delta }\left( t\right) \left( x,y\right) =S\left(
t,A\right) U_{\Delta }\left( t\right) \left( x,y\right) 
\]%
where $C\left( t,A\right) $ is a cosine, $S\left( t,A\right) $ is a sine
operator-functions (see e.g. $\left[ 11\right] $) with generator of $A$, i.e.%
\[
C\left( t,A\right) =\frac{1}{2}\left( e^{tA^{\frac{1}{2}}}+e^{-tA^{\frac{1}{2%
}}}\right) \text{, }S\left( t,A\right) =\frac{1}{2}A^{-\frac{1}{2}}\left(
e^{tA^{\frac{1}{2}}}-e^{-tA^{\frac{1}{2}}}\right) . 
\]%
and $U_{\Delta }\left( t\right) \left( x,y\right) =e^{\Delta t}\left(
x,y\right) $ is a fundamental solution of the free wave equation:%
\[
\partial _{t}^{2}u-\Delta u=0,\text{ }x\in R^{n},\text{ }t\in \left[ 0,T%
\right] , 
\]%
i.e. 
\begin{equation}
U_{\Delta }\left( t\right) \left( x,y\right) =\left( 4\pi t\right) ^{-\frac{n%
}{2}}e^{\left\vert x-y\right\vert ^{2}\mid 4t}\text{, }t\neq 0,  \tag{3.3}
\end{equation}

\[
U_{\Delta }\left( t\right) f\left( x\right) =\left( 2\pi t\right) ^{-\frac{n%
}{2}}\dint\limits_{R^{n}}e^{\frac{\left\vert x-y\right\vert ^{2}}{2t}%
}f\left( y\right) dy. 
\]

\textbf{Lemma 3.1.\ }Let $A$ be an absolute positive operator in a Banach\
space $E$ and $0\leq \alpha <1$. Then the following dispersive inequalites
hold 
\begin{equation}
\left\Vert A^{\alpha }U_{\Delta +A}\left( t\right) f\right\Vert
_{L_{x}^{p}\left( R^{n}:E\right) }\lesssim t^{-\left[ n\left( \frac{1}{2}-%
\frac{1}{p}\right) +\alpha \right] }\left\Vert f\right\Vert
_{L_{x}^{p^{\prime }}\left( R^{n}:E\right) },  \tag{3.4}
\end{equation}%
\begin{equation}
\left\Vert A^{\alpha }U_{\Delta +A}\left( t-s\right) f\right\Vert
_{L^{\infty }\left( R^{n};E\right) }\lesssim \left\vert t-s\right\vert
^{-\left( \frac{n}{2}+\alpha \right) }\left\Vert f\right\Vert _{L^{1}\left(
R^{n};E\right) }  \tag{3.5}
\end{equation}%
for $t\neq 0,$ $2\leq p\leq \infty ,$ $\frac{1}{p}+\frac{1}{p^{\prime }}=1.$

\textbf{Proof. }By using $\left( 3.3\right) $ and Young's integral
inequality we have 
\begin{equation}
\left\Vert U_{\Delta }\left( t\right) f\right\Vert _{L_{x}^{p}\left(
R^{n}:E\right) }\lesssim \left\vert t\right\vert ^{-n\left( \frac{1}{2}-%
\frac{1}{p}\right) }\left\Vert f\right\Vert _{L_{x}^{p^{\prime }}\left(
R^{n}:E\right) },  \tag{3.6}
\end{equation}%
\[
\left\Vert U_{\Delta }\left( t\right) f\right\Vert _{L_{x}^{\infty }\left(
R^{n}:E\right) }\lesssim \left\vert t\right\vert ^{-\frac{n}{2}}\left\Vert
f\right\Vert _{L_{x}^{1}\left( R^{n}:E\right) }. 
\]%
By $\left( 1.10\right) $ we get 
\[
\left\Vert A^{\alpha }U_{A}\left( t\right) \right\Vert _{B\left( E\right)
}\lesssim \left\vert t\right\vert ^{-\alpha }\text{, }t\neq 0. 
\]%
By using then the properties of $U_{\Delta +A}\left( t\right) =U_{\Delta
}\left( t\right) $ $U_{A}\left( t\right) $, the estimates $\left( 3.7\right) 
$ and $\left( 3.6\right) $ we obtain $\left( 3.4\right) $ and $\left(
3.5\right) .$

In this section, we make the necessary estimates to solution of the
following Cauchy problem for the linear abstract wave equation%
\begin{equation}
u_{tt}-\Delta u+Au=F\left( t,x\right) ,\text{ }x\in R^{n},\text{ }t\in
\left( 0,\infty \right) ,  \tag{3.7}
\end{equation}%
\begin{equation}
u\left( 0,x\right) =\varphi \left( x\right) +\dsum\limits_{k=1}^{m}\alpha
_{k}u\left( \lambda _{k},x\right) ,\text{ for a.e. }x\in R^{n},  \tag{3.8}
\end{equation}%
\begin{equation}
u_{t}\left( 0,x\right) =\psi \left( x\right) +\dsum\limits_{k=1}^{m}\beta
_{k}u_{t}\left( \lambda _{k},x\right) ,\text{ for a.e. }x\in R^{n}, 
\tag{3.9}
\end{equation}%
where $A$ is a linear operator in a Hilbert space $H.$

\textbf{Condition 3.1. }Assume:

(1) $\left\vert \alpha _{k}+\beta _{k}\right\vert >0$, $\dsum%
\limits_{k,j=1}^{m}\alpha _{k}\beta _{j}\neq 0;$

(2) $E$ is a Banach space;

(3) $A$ is absulute positive operator in a Banach\ space $E$ and $\gamma
\geq \frac{n}{p}$ for $p\in \left[ 1,\infty \right] .$

First we need the following lemma:

\textbf{Lemma 3.2. }Suppose the Condition 3.1 hold, $\varphi \in \mathring{W}%
^{\gamma ,p}\left( R^{n};E\left( A\right) \right) $ and $\psi \in \mathring{W%
}^{\gamma -1,p}\left( R^{n};E\left( A\right) \right) $. Then problem $\left(
2.1\right) -\left( 2.2\right) $ has a unique generalized solution.

\textbf{Proof. }By using of the Fourier transform we get from $(3.1):$%
\[
\hat{u}_{tt}\left( t,\xi \right) +A_{\xi }\hat{u}\left( t,\xi \right) =\hat{F%
}\left( t,\xi \right) ,\text{ } 
\]%
\begin{equation}
\hat{u}\left( 0,x\right) =\hat{\varphi}\left( \xi \right)
+\dsum\limits_{k=1}^{m}\alpha _{k}\hat{u}\left( \lambda _{k},\xi \right) ,%
\text{ }  \tag{3.10}
\end{equation}%
\[
\hat{u}_{t}\left( 0,\xi \right) =\hat{\psi}\left( \xi \right)
+\dsum\limits_{k=1}^{m}\beta _{k}\hat{u}_{t}\left( \lambda _{k},\xi \right) 
\text{ for a.e. }\xi \in R^{n}. 
\]%
where $\hat{u}\left( t,\xi \right) $ is a Fourier transform of $u\left(
t,x\right) $ with respect to $x$ and%
\[
A_{\xi }=A+\left\vert \xi \right\vert ^{2}\text{, }\xi \in R^{n}. 
\]%
Consider the problem%
\begin{equation}
\hat{u}_{tt}\left( t,\xi \right) +A_{\xi }\hat{u}\left( t,\xi \right) =\hat{F%
}\left( t,\xi \right) ,\text{ }  \tag{3.11}
\end{equation}%
\[
\hat{u}\left( 0,\xi \right) =u_{0}\left( \xi \right) ,\text{ }\hat{u}%
_{t}\left( 0,\xi \right) =u_{1}\left( \xi \right) ,\text{ }\xi \in R^{n},%
\text{ }t\in \left[ 0,T\right] ,\text{ } 
\]%
where $u_{0}\left( \xi \right) $, $u_{1}\left( \xi \right) \in D\left(
A\right) $ and for $\xi \in R^{n}.$ By virtue of $\left[ \text{11, \S 11.2,
11.4}\right] $ we obtain that $A_{\xi }$ is a generator of a strongly
continuous cosine operator function\ and the Cauchy problem $(3.11)$ has a
unique solution for all $\xi \in R^{n},$ moreover, the solution can be
expressed as%
\begin{equation}
\hat{u}\left( t,\xi \right) =C\left( t,\xi ,A\right) u_{0}\left( \xi \right)
+S\left( t,\xi ,A\right) u_{1}\left( \xi \right) +  \tag{3.12}
\end{equation}%
\[
\dint\limits_{0}^{t}S\left( t-\tau ,\xi ,A\right) \hat{F}\left( \tau ,\xi
\right) d\tau ,\text{ }t\in \left( 0,T\right) , 
\]%
where $C\left( t,\xi ,A\right) $ is a cosine and $S\left( t,\xi ,A\right) $
is a sine operator-functions generated by $A_{\xi }$, i.e.%
\[
C\left( t,\xi ,A\right) =\frac{1}{2}\left( e^{tA_{\xi }^{\frac{1}{2}%
}}+e^{-tA_{\xi }^{\frac{1}{2}}}\right) \text{, }S\left( t,\xi ,A\right) =%
\frac{1}{2}A_{\xi }^{-\frac{1}{2}}\left( e^{tA_{\xi }^{\frac{1}{2}%
}}-e^{-tA_{\xi }^{\frac{1}{2}}}\right) . 
\]%
Using the formula $\left( 3.12\right) $ and the condition $\left(
3.10\right) $ we get 
\[
u_{0}\left( \xi \right) =\hat{\varphi}\left( \xi \right)
+\dsum\limits_{k=1}^{m}\alpha _{k}\left[ C\left( \lambda _{k},\xi ,A\right)
u_{0}\left( \xi \right) +S\left( \lambda _{k},\xi ,A\right) u_{1}\left( \xi
\right) \right] + 
\]

\[
\dsum\limits_{k=1}^{m}\alpha _{k}\dint\limits_{0}^{\lambda _{k}}S\left(
\lambda _{k}-\tau ,\xi ,A\right) \hat{F}\left( \tau ,\xi \right) d\tau ,%
\text{ }\tau \in \left( 0,T\right) . 
\]%
Then, 
\[
\left[ I-\dsum\limits_{k=1}^{m}\alpha _{k}C\left( \lambda _{k},\xi ,A\right) %
\right] u_{0}\left( \xi \right) -\dsum\limits_{k=1}^{m}\alpha _{k}S\left(
\lambda _{k},\xi ,A\right) u_{1}\left( \xi \right) = 
\]%
\begin{equation}
\dsum\limits_{k=1}^{m}\alpha _{k}\dint\limits_{0}^{\lambda _{k}}S\left(
\lambda _{k}-\tau ,\xi ,A\right) \hat{F}\left( \tau ,\xi \right) d\tau +\hat{%
\varphi}\left( \xi \right) .  \tag{3.13}
\end{equation}%
Differentiating both sides of formula $\left( 3.12\right) $ we obtain%
\[
\hat{u}_{t}\left( t,\xi \right) =A_{\xi }S\left( t,\xi ,A\right) u_{0}\left(
\xi \right) +C\left( t,\xi ,A\right) u_{1}\left( \xi \right) + 
\]%
\[
\frac{1}{2}\hat{F}\left( t,\xi \right) +\dint\limits_{0}^{t}AS\left( t-\tau
,\xi ,A\right) \hat{F}\left( \tau ,\xi \right) d\tau ,\text{ }t\in \left(
0,\infty \right) . 
\]%
Using the above formula and the condition 
\[
\hat{u}_{t}\left( 0,\xi \right) =\hat{\psi}\left( \xi \right)
+\dsum\limits_{k=1}^{m}\beta _{k}\hat{u}_{t}\left( \lambda _{k},\xi \right) 
\]%
we obtain%
\[
u_{1}\left( \xi \right) =\hat{\psi}\left( \xi \right)
+\dsum\limits_{k=1}^{m}\beta _{k}\left[ A_{\xi }S\left( \lambda _{k},\xi
,A\right) u_{0}\left( \xi \right) +C\left( \lambda _{k},\xi ,A\right)
u_{1}\left( \xi \right) \right] + 
\]

\[
\dsum\limits_{k=1}^{m}\beta _{k}\left[ \frac{1}{2}\hat{g}\left( \lambda
_{k},\xi \right) +\dint\limits_{0}^{\lambda _{k}}S\left( \lambda _{k}-\tau
,\xi ,A\right) \hat{F}\left( \tau ,\xi \right) d\tau \right] . 
\]%
Thus, 
\[
-\dsum\limits_{k=1}^{m}\beta _{k}A_{\xi }S\left( \lambda _{k},\xi ,A\right)
u_{0}\left( \xi \right) +\left[ I-\dsum\limits_{k=1}^{m}\beta _{k}C\left(
\lambda _{k},\xi ,A\right) \right] u_{1}\left( \xi \right) = 
\]%
\begin{equation}
\dsum\limits_{k=1}^{m}\beta _{k}\left[ \frac{1}{2}\hat{g}\left( \lambda
_{k},\xi \right) +\dint\limits_{0}^{\lambda _{k}}S\left( \lambda _{k}-\tau
,\xi ,A\right) \hat{F}\left( \tau ,\xi \right) d\tau \right] +\hat{\psi}%
\left( \xi \right) .  \tag{3.14}
\end{equation}%
Now, we consider the system of equations $\left( 3.13\right) $ and $\left(
3.14\right) $ in $u_{0}\left( \xi \right) $ and $u_{1}\left( \xi \right) $.
The determinant of this system is 
\[
D\left( \xi \right) =\left\vert 
\begin{array}{cc}
\alpha _{11}\left( \xi \right) & \alpha _{12}\left( \xi \right) \\ 
\alpha _{21}\left( \xi \right) & \alpha _{22}\left( \xi \right)%
\end{array}%
\right\vert , 
\]%
where 
\[
\alpha _{11}\left( \xi \right) =I-\dsum\limits_{k=1}^{m}\alpha _{k}C\left(
\lambda _{k},\xi ,A\right) ,\text{ }\alpha _{12}\left( \xi \right)
=-\dsum\limits_{k=1}^{m}\alpha _{k}S\left( \lambda _{k},\xi ,A\right) , 
\]%
\[
\alpha _{21}\left( \xi \right) =-\dsum\limits_{k=1}^{m}\beta _{k}A_{\xi
}S\left( \lambda _{k},\xi ,A\right) ,\text{ }\alpha _{22}\left( \xi \right)
=I-\dsum\limits_{k=1}^{m}\beta _{k}C\left( \lambda _{k},\xi ,A\right) . 
\]%
We find the determinant of the system $\left( 3.13\right) $-$\left(
3.14\right) :$%
\[
D\left( \xi \right) =I-\dsum\limits_{k=1}^{m}\left( \alpha _{k}+\beta
_{k}\right) C\left( \lambda _{k},\xi ,A\right) + 
\]%
\[
\dsum\limits_{k,j=1}^{m}\alpha _{k}\beta _{j}\left[ C\left( \lambda _{k},\xi
,A\right) C\left( \lambda _{j},\xi ,A\right) -A_{\xi }S\left( \lambda
_{k},\xi ,A\right) S\left( \lambda _{j},\xi ,A\right) \right] . 
\]

By properties of operator functions $C\left( \lambda ,\xi ,A\right) $ and $%
S\left( \lambda ,\xi ,A\right) $ we get $D\left( \xi \right) \neq 0.$
Solving the system $\left( 3.13\right) -\left( 3.14\right) $, by using the
property of sine and cosine operator function $\left[ \text{11, \S 11.2, 11.4%
}\right] $\ we get%
\begin{equation}
u_{0}\left( \xi \right) =D^{-1}\left( \xi \right) \left[ \left(
I-\dsum\limits_{k=1}^{m}\beta _{k}C\left( \lambda _{k},\xi ,A\right) \right)
f_{1}-\dsum\limits_{k=1}^{m}\alpha _{k}S\left( \lambda _{k},\xi ,A\right)
f_{2}\right] ,  \tag{3.15}
\end{equation}%
\[
u_{1}\left( \xi \right) =D^{-1}\left( \xi \right) \left[ \left(
I-\dsum\limits_{k=1}^{m}\alpha _{k}C\left( \lambda _{k},\xi ,A\right)
\right) f_{2}+\dsum\limits_{k=1}^{m}\beta _{k}A_{\xi }S\left( \lambda
_{k},\xi ,A\right) f_{1}\right] , 
\]%
where 
\[
f_{1}=\dsum\limits_{k=1}^{m}\alpha _{k}\dint\limits_{0}^{\lambda
_{k}}S\left( \lambda _{k}-\tau ,\xi ,A\right) \hat{F}\left( \tau ,\xi
\right) d\tau +\hat{\varphi}\left( \xi \right) , 
\]%
\begin{equation}
f_{2}=\dsum\limits_{k=1}^{m}\beta _{k}\left[ \frac{1}{2}\hat{F}\left(
\lambda _{k},\xi \right) +\dint\limits_{0}^{\lambda _{k}}S\left( \lambda
_{k}-\tau ,\xi ,A\right) \hat{F}\left( \tau ,\xi \right) d\tau \right] +\hat{%
\psi}\left( \xi \right) .  \tag{3.16}
\end{equation}%
From $\left( 3.12\right) ,$ $\left( 3.15\right) $ and $\left( 3.16\right) $
we get that the solution of $\left( 3.10\right) $ can be expressed as 
\[
\hat{u}\left( t,\xi \right) =D^{-1}\left( \xi \right) \left\{ \left[ C\left(
t,\xi ,A\right) \left( I-\dsum\limits_{k=1}^{m}\beta _{k}C\left( \lambda
_{k},\xi ,A\right) \right) \right. \right. + 
\]%
\[
\left. S\left( t,\xi ,A\right) \dsum\limits_{k=1}^{m}\beta _{k}A_{\xi
}S\left( \lambda _{k},\xi ,A\right) \right] f_{1}+\left[ S\left( t,\xi
,A\right) \left( I-\dsum\limits_{k=1}^{m}\alpha _{k}C\left( \lambda _{k},\xi
,A\right) \right) \right. - 
\]%
\begin{equation}
\left. \left. C\left( t,\xi ,A\right) \dsum\limits_{k=1}^{m}\alpha
_{k}S\left( \lambda _{k},\xi ,A\right) \right] f_{2}\right\} ,\text{ }t\in
\left( 0,T\right) .  \tag{3.17}
\end{equation}%
We obtain from $\left( 3.17\right) $ that there is a generalized solution of 
$(3.7)-(3.9)$ given by 
\begin{equation}
u\left( t,x\right) =S_{1}\left( t,x,A\right) \varphi \left( x\right)
+S_{2}\left( t,x,A\right) \psi \left( x\right) +\Phi \left( t,x,A\right) , 
\tag{3.18}
\end{equation}%
where $S_{1}\left( t,A\right) $ and $S_{2}\left( t,A\right) $ are linear
operator functions in $E$ defined by 
\[
S_{1}\left( t,x,A\right) \varphi =F^{-1}D^{-1}\left( \xi \right) \left[
C\left( t,\xi ,A\right) \left( I-\dsum\limits_{k=1}^{m}\beta _{k}C\left(
\lambda _{k},\xi ,A\right) \right) \right. +\text{ } 
\]%
\[
\left. S\left( t,\xi ,A\right) \dsum\limits_{k=1}^{m}\beta _{k}A_{\xi
}S\left( \lambda _{k},\xi ,A\right) \right] \hat{\varphi}\left( \xi \right)
, 
\]

\begin{equation}
S_{2}\left( t,x,A\right) \psi =F^{-1}D^{-1}\left( \xi \right) \left[
-C\left( t,\xi ,A\right) \dsum\limits_{k=1}^{m}\alpha _{k}S\left( \lambda
_{k},\xi ,A\right) \right. +  \tag{3.19}
\end{equation}%
\[
\left. S\left( t,\xi ,A\right) \left( I-\dsum\limits_{k=1}^{m}\alpha
_{k}C\left( \lambda _{k},\xi ,A\right) \right) \right] \hat{\psi}\left( \xi
\right) , 
\]

\[
\Phi \left( t,x,A\right) =F^{-1}D^{-1}\left( \xi \right) \left\{ \left[
C\left( t,\xi ,A\right) \right. \left( I-\dsum\limits_{k=1}^{m}\beta
_{k}C\left( \lambda _{k},\xi ,A\right) \right) g_{1}\left( t\right) \right.
- 
\]

\[
\left. \dsum\limits_{k=1}^{m}\alpha _{k}S\left( \lambda _{k},\xi ,A\right)
g_{2}\left( \xi \right) \right] +S\left( t,\xi ,A\right) \left[
\dsum\limits_{k=1}^{m}\beta _{k}A_{\xi }S\left( \lambda _{k},\xi ,A\right)
g_{1}\left( \xi \right) \right. + 
\]%
\[
\left( I-\dsum\limits_{k=1}^{m}\alpha _{k}C\left( \lambda _{k},\xi ,A\right)
\right) \left. \dsum\limits_{k=1}^{m}\beta _{k}A_{\xi }S\left( \lambda
_{k},\xi ,A\right) g_{2}\left( \xi \right) \right\} , 
\]%
here 
\begin{equation}
g_{1}\left( \xi \right) =\dsum\limits_{k=1}^{m}\alpha
_{k}\dint\limits_{0}^{\lambda _{k}}S\left( \lambda _{k}-\tau ,\xi ,A\right) 
\hat{F}\left( \tau ,\xi \right) d\tau ,  \tag{3.20}
\end{equation}%
\[
g_{2}\left( \xi \right) =\dsum\limits_{k=1}^{m}\beta _{k}\left[ \frac{1}{2}%
\hat{F}\left( \lambda _{k},\xi \right) +\dint\limits_{0}^{\lambda
_{k}}S\left( \lambda _{k}-\tau ,\xi ,A\right) \hat{F}\left( \tau ,\xi
\right) d\tau \right] . 
\]

\begin{center}
\bigskip \textbf{4}. \textbf{Strichartz inequalities} \textbf{for linear
wave equation}
\end{center}

\bigskip The proof of Strichartz type estimates involves basically two type
ingredients. The first one consists of specific estimates, in particular
stationary phase estimates, on evolution groups associated with homogenous
equations. The second one consists of abstract arguments, not specific to
the wave equations. This is mainly duality argument and were first applied
in $\left[ 37\right] .$

\textbf{Condition 4.1. }Assume $n>1,$%
\[
\frac{1}{q}+\frac{n-1}{2r}\leq \frac{n-1}{4},\text{ }2\leq q,r\leq \infty 
\text{ \ and }\left( n,\text{ }q,\text{ }r\right) \neq \left( 2,\text{ }2,%
\text{ }\infty \right) . 
\]

\textbf{Remark 4.1. }If $\frac{1}{q}+\frac{n-1}{2r}=\frac{n-1}{4},$ then $(q,
$ $r)$ is called sharp admissible, otherwise $(q,$ $r)$ is called nonsharp
admissible. Note in particular that when $n>2$ the endpoint $\left( 2\text{, 
}\frac{2\left( n-1\right) }{n-3}\right) $ is called sharp admissible.

\bigskip For a space-time slab $\left[ 0,T\right] \times R^{n}$, we define
the $E-$valued Strichartz norm%
\[
\left\Vert u\right\Vert _{S^{0}\left( \left[ 0,T\right] ;E\right)
}=\sup\limits_{\left( q,r\right) \text{ admissible}}\left\Vert u\right\Vert
_{L_{t}^{q}L_{x}^{r}\left( \left[ 0,T\right] \times R^{n};E\right) }, 
\]%
where $S^{0}\left( \left[ 0,T\right] ;E\right) $ is the closure of all $E-$%
valued test functions under this norm and $N^{0}\left( \left[ 0,T\right]
;E\right) $ denotes the dual of $S^{0}\left( \left[ 0,T\right] ;E\right) .$

Assume $H$ is an abstract Hilbert space and $Q$ is a $H-$valued Hilbert
space of function. Suppose for each $t\in \mathbb{R}$ an operator $U\left(
t\right) $: $Q\rightarrow L^{2}\left( \Omega ;E\right) $ obeys the following
estimates:

\begin{equation}
\left\Vert U\left( t\right) f\right\Vert _{L_{x}^{2}\left( \Omega ;H\right)
}\lesssim \left\Vert f\right\Vert _{Q}  \tag{4.1}
\end{equation}%
for all $t,$ $\Omega \subset R^{n}$ and all $f\in Q.$ Moreover,%
\begin{equation}
\left\Vert U\left( s\right) U^{\ast }\left( t\right) g\right\Vert
_{L_{x}^{\infty }\left( \Omega ;H\right) }\lesssim \left\vert t-s\right\vert
^{-\frac{n-1}{2}}\left\Vert g\right\Vert _{L_{x}^{1}\left( \Omega ;H\right) }
\tag{4.2}
\end{equation}%
\begin{equation}
\left\Vert U\left( s\right) U^{\ast }\left( t\right) g\right\Vert
_{L_{x}^{\infty }\left( \Omega ;H\right) }\lesssim \left( 1+\left\vert
t-s\right\vert ^{-\frac{n-1}{2}}\right) \left\Vert g\right\Vert
_{L_{x}^{1}\left( \Omega ;H\right) }  \tag{4.3}
\end{equation}%
for all $t\neq s$ and all $g\in L_{x}^{1}\left( \Omega ;H\right) .$

For proving the main theorem of this section, we will use the following
bilinear interpolation result (see $\left[ 5\right] $, Section 3.13.5(b)).

\bigskip \textbf{Lemma 4.1. }Assume $A_{0}$, $A_{1},$ $B_{0}$, $B_{1},$ $%
C_{0}$, $C_{1}$ are Banach spaces and $T$ is a bilinear operator bounded
from ($A_{0}\times B_{0}$, $A_{0}\times B_{1},$ $A_{1}\times B_{0}$ ) into ($%
C_{0}$, $C_{1}$, $C_{1}$), respectively. Then whenever $0<\theta _{0},$ $%
\theta _{1}<\theta <1$ are such that $1\leq \frac{1}{p}+\frac{1}{q}$ and $%
\theta =\theta _{0}+$ $\theta _{1}$, the operator is bounded from 
\[
\left( A_{0}\text{, }A_{1}\right) _{\theta _{0}pr}\times \left( B_{0}\text{, 
}B_{1}\right) _{\theta _{1}qr} 
\]%
to $\left( C_{0}\text{, }C_{1}\right) _{\theta r}.$

By following $\left[ \text{22, Theorem 1.2}\right] $ we have:

\textbf{Theorem 4.1. }Assume $U(t)$ obeys $\left( 4.1\right) $-$\left(
4.3\right) $. Let $U\left( t\right) $ generates absolute positive
infinitesimal generator operator $A$ and $0\leq \alpha <1.$ Then the
following estimates are hold%
\begin{equation}
\left\Vert U\left( t\right) f\right\Vert _{L_{t}^{q}L_{x}^{r}\left( H\right)
}\lesssim \left\Vert f\right\Vert _{Q},  \tag{4.4}
\end{equation}

\begin{equation}
\left\Vert \dint U^{\ast }\left( s\right) F\left( s\right) ds\right\Vert
_{Q}\lesssim \left\Vert F\right\Vert _{L_{t}^{q^{\prime }}L_{x}^{r^{\prime
}}\left( H\right) },  \tag{4.5}
\end{equation}%
\begin{equation}
\dint\limits_{s<t}\left\Vert A^{\alpha }U\left( t\right) U^{\ast }\left(
s\right) F\left( s\right) ds\right\Vert _{L_{t}^{q}L_{x}^{r}\left( H\right)
}\lesssim \left\Vert F\right\Vert _{L_{t}^{\tilde{q}^{\prime }}L_{x}^{\tilde{%
r}^{\prime }}\left( H\right) },  \tag{4.6}
\end{equation}%
for all sharp admissible exponent pairs $\left( q,r\right) $, $\left( \tilde{%
q},\tilde{r}\right) .$ Furthermore, if the decay hypothesis is strengthened
to $(4.3)$, then $(4.4)-\left( 4.6\right) $ hold for all admissible $\left(
q,\text{ }r\right) $, $\left( \tilde{q},\text{ }\tilde{r}\right) .$

\textbf{Proof. The first step: }Consider the nonendpoint case, i.e. $\left(
q,\text{ }r\right) \neq $ $\left( 2,\text{ }\frac{2\left( n-1\right) }{n-3}%
\right) $ and will show firstly, the estimates $\left( 4.4\right) $, $\left(
4.5\right) .$ By duality, $(4.4)$ is equivalent to $(4.5)$. By the $TT^{\ast
}$ method, $(4.6)$ is in turn equivalent to the bilinear form estimate 
\begin{equation}
\left\vert \dint \dint \langle \left( A^{\frac{\alpha }{2}}U\left( s\right)
\right) ^{\ast }F\left( s\right) ,\left( A^{\frac{\alpha }{2}}U\left(
t\right) \right) ^{\ast }G\left( t\right) \rangle dsdt\right\vert \lesssim 
\tag{4.7}
\end{equation}

\[
\left\Vert F\right\Vert _{L_{t}^{q^{\prime }}L_{x}^{r^{\prime }}\left(
H\right) }\left\Vert G\right\Vert _{L_{t}^{q^{\prime }}L_{x}^{r^{\prime
}}\left( H\right) }. 
\]%
By symmetry it suffices to show to the retarded version of $\left(
4.7\right) $%
\begin{equation}
\left\vert T\left( F,G\right) \right\vert \lesssim \left\Vert F\right\Vert
_{L_{t}^{q^{\prime }}L_{x}^{r^{\prime }}\left( H\right) }\left\Vert
G\right\Vert _{L_{t}^{q^{\prime }}L_{x}^{r^{\prime }}\left( H\right) }, 
\tag{4.8}
\end{equation}%
where $T\left( F,G\right) $ is the bilinear form defined by 
\[
T\left( F,G\right) =\dint \dint\limits_{s<t}\langle \left( A^{\frac{\alpha }{%
2}}U\left( s\right) \right) ^{\ast }F\left( s\right) ,\left( A^{\frac{\alpha 
}{2}}U\left( t\right) \right) ^{\ast }G\left( t\right) \rangle dsdt 
\]

By real interpolation between the bilinear form of $\left( 4.1\right) $ and
due to estimate $\left( 1.10\right) $ we get 
\[
\left\vert \langle \left( A^{\frac{\alpha }{2}}U\left( s\right) \right)
^{\ast }F\left( s\right) ,\left( A^{\frac{\alpha }{2}}U\left( t\right)
\right) ^{\ast }G\left( t\right) \rangle \right\vert \lesssim \left\Vert
F\left( s\right) \right\Vert _{L_{x}^{2}}\left\Vert G\left( t\right)
\right\Vert _{L_{x}^{2}}. 
\]%
By using the bilinear form of $\left( 4.2\right) $ and $\left( 1.10\right) $
we have 
\begin{equation}
\left\vert \langle \left( A^{\frac{\alpha }{2}}U\left( s\right) \right)
^{\ast }F\left( s\right) ,\left( A^{\frac{\alpha }{2}}U\left( t\right)
\right) ^{\ast }G\left( t\right) \rangle \right\vert \lesssim  \tag{4.9}
\end{equation}%
\[
\left\vert t-s\right\vert ^{-\frac{n}{2}}\left\Vert F\left( s\right)
\right\Vert _{L_{x}^{1}\left( \Omega ;H\right) }\left\Vert G\left( t\right)
\right\Vert _{L_{x}^{1}\left( \Omega ;H\right) }. 
\]%
In a similar way, we obtain 
\begin{equation}
\left\vert \langle \left( A^{\frac{\alpha }{2}}U\left( s\right) \right)
^{\ast }F\left( s\right) ,\left( A^{\frac{\alpha }{2}}U\left( t\right)
\right) ^{\ast }G\left( t\right) \rangle \right\vert \lesssim  \tag{4.10}
\end{equation}%
\[
\left\vert t-s\right\vert ^{--1-\beta \left( r,r\right) }\left\Vert F\left(
s\right) \right\Vert _{L_{x}^{r^{\prime }}\left( \Omega ;H\right)
}\left\Vert G\left( t\right) \right\Vert _{L_{x}^{r^{\prime }}\left( \Omega
;H\right) }, 
\]%
where $\beta (r,\tilde{r})$ is given by 
\begin{equation}
\beta (r,\tilde{r})=\frac{n}{2}-1-\frac{n}{2}\left( \frac{1}{r}-\frac{1}{%
\tilde{r}}\right) .  \tag{4.11}
\end{equation}

It is clear that $\beta (r,r)\leq 0$ when $n>2.$ In the sharp admissible
case we have 
\[
\frac{1}{q}+\frac{1}{q^{\prime }}=-\beta (r,r), 
\]%
and $\left( 4.8\right) $ follows from $\left( 4.10\right) $ and the
Hardy-Littlewood-Sobolev inequality ($[20]$) when $q>q^{\prime }.$

If we are assuming the truncated decay $(4.3)$, then $(4.10)$ can be
improved to 
\begin{equation}
\left\vert \langle \left( A^{\frac{\alpha }{2}}U\left( s\right) \right)
^{\ast }F\left( s\right) ,\left( A^{\frac{\alpha }{2}}U\left( t\right)
\right) ^{\ast }G\left( t\right) \rangle \right\vert \lesssim  \tag{4.12}
\end{equation}%
\[
\left( 1+\left\vert t-s\right\vert \right) ^{-1-\beta \left( r,r\right)
}\left\Vert F\left( s\right) \right\Vert _{L_{x}^{r^{\prime }}\left( \Omega
;H\right) }\left\Vert G\left( t\right) \right\Vert _{L_{x}^{r^{\prime
}}\left( \Omega ;H\right) } 
\]%
and now Young's inequality gives $(4.8)$ when 
\[
-\beta (r,r)+\frac{1}{q}>\frac{1}{q^{\prime }}, 
\]%
i.e. $(q,r)$ is nonsharp admissible. This concludes the proof of $\left(
4.4\right) $ and $(4.5)$ for nonendpoint case.

\textbf{The second step; }It remains to prove $\left( 4.4\right) $ and $(4.5)
$ for the endpoint case, i.e. when%
\[
\left( q,\text{ }r\right) =\left( 2,\text{ }\frac{2\left( n-1\right) }{n-3}%
\right) ,\text{ }n>2.
\]%
It suffices to show $(4.8)$. \ By decomposing $T(F,G)$ dyadically as $%
\dsum\limits_{j}T_{j}(F,G),$ where the summation is over the integers $%
\mathbb{Z}$ and%
\begin{equation}
T_{j}\left( F,G\right) =\dint\limits_{t-2^{j-1}<s\leq t-2^{j}}\langle \left(
A^{\frac{\alpha }{2}}U\left( s\right) \right) ^{\ast }F\left( s\right)
,\left( A^{\frac{\alpha }{2}}U\left( t\right) \right) ^{\ast }G\left(
t\right) \rangle dsdt  \tag{4.13}
\end{equation}%
we see that it suffices to prove the estimate 
\begin{equation}
\dsum\limits_{j}\left\vert T_{j}(F,G)\right\vert \lesssim \left\Vert
F\right\Vert _{L_{t}^{2}L_{x}^{r^{\prime }}\left( H\right) }\left\Vert
G\right\Vert _{L_{t}^{2}L_{x}^{r^{\prime }}\left( H\right) }\text{.} 
\tag{4.14}
\end{equation}%
For this aim, before we will show the following estimate 
\begin{equation}
\left\vert T_{j}(F,G)\right\vert \lesssim 2^{-j\beta \left( a,b\right)
}\left\Vert F\right\Vert _{L_{t}^{2}L_{x}^{a^{\prime }}\left( H\right)
}\left\Vert G\right\Vert _{L_{t}^{2}L_{x}^{b^{\prime }}\left( H\right) } 
\tag{4.15}
\end{equation}%
for all $j\in \mathbb{Z}$ and all $\left( \frac{1}{a},\frac{1}{b}\right) $\
in a neighbourhood of $\left( \frac{1}{r},\frac{1}{r}\right) $. For proving $%
\left( 4.15\right) $ we will use the real interpolation of $H$-valued
Lebesque space and sequence spaces $l_{q}^{s}\left( H\right) $ (see e.g $%
\left[ \text{38}\right] $ \S\ 1.18.2 and 1.18.6). Indeed, by $\left[ \text{%
38, \S\ 1.18.4.}\right] $ we have 
\begin{equation}
\left( L_{t}^{2}L_{x}^{p_{0}}\left( H\right) ,L_{t}^{2}L_{x}^{p_{1}}\left(
H\right) \right) _{\theta ,2}=L_{t}^{2}L_{x}^{p,2}\left( H\right)  
\tag{4.16}
\end{equation}%
whenever $p_{0},$ $p_{1}\in \left[ 1,\infty \right] ,$ $p_{0}\neq p_{1}$ and 
$\frac{1}{p}=\frac{1-\theta }{p_{0}}+\frac{\theta }{p_{1}}$ and $\left(
l_{\infty }^{s_{0}}\left( H\right) ,l_{\infty }^{s_{1}}\left( H\right)
\right) _{\theta ,1}=l_{1}^{s}\left( H\right) $ for $s_{0}$, $s_{1}\in 
\mathbb{R}$, $s_{0}\neq s_{1}$ and%
\[
\frac{1}{s}=\frac{1-\theta }{s_{0}}+\frac{\theta }{s_{1}},
\]%
where 
\[
l_{q}^{s}\left( H\right) =\left\{ u=\left\{ u_{j}\right\} _{j=1}^{\infty
},u_{j}\in H\text{,}\right. \text{ }
\]%
\[
\left\Vert u\right\Vert _{l_{q}^{s}\left( H\right) }=\left. \left(
\dsum\limits_{j=1}^{\infty }2^{jsq}\left\Vert u_{j}\right\Vert
_{H}^{q}\right) ^{\frac{1}{q}}<\infty \right\} .
\]

By $\left( 4.16\right) $ the estimate $(4.15)$ can be rewritten as 
\begin{equation}
T:L_{t}^{2}L_{x}^{a^{\prime }}\left( H\right) \times
L_{t}^{2}L_{x}^{b^{\prime }}\left( H\right) \rightarrow l_{\infty }^{\beta
\left( a,b\right) },  \tag{4.17}
\end{equation}%
where $T=\left\{ T_{j}\right\} $ is the vector-valued bilinear operator
corresponding to the $T_{j}.$ We apply Lemma 3.2 to $\left( 4.17\right) $
with $r=1$, $p=q=2$ and arbitrary exponents $a_{0},$ $a_{1}$, $b_{0}$, $%
b_{1} $ such that 
\[
\beta \left( a_{0},b_{1}\right) =\beta \left( a_{1},b_{0}\right) \neq \beta
\left( a_{0},b_{0}\right) . 
\]

Using the real interpolation space identities we obtain%
\[
T:L_{t}^{2}L_{x}^{a^{\prime },2}\left( E^{\ast }\right) \times
L_{t}^{2}L_{x}^{b^{\prime },2}\left( E^{\ast }\right) \rightarrow
l_{1}^{\beta \left( a,b\right) } 
\]%
for all $(a,b)$ in a neighbourhood of $(r,r)$. Applying this to $a=b=r$ and
using the fact that $L^{r^{\prime }}\left( H\right) \subset L^{r^{\prime
},2}\left( H\right) $ we obtain%
\[
T:L_{t}^{2}L_{x}^{a^{\prime },2}\left( H\right) \times
L_{t}^{2}L_{x}^{b^{\prime },2}\left( H\right) \rightarrow l_{1}^{0}\left(
H\right) 
\]%
which implies $\left( 4.15\right) .$

We are now ready to state the Strichartz estimates to solution $\left(
3.7\right) -\left( 3.9\right) $.

\textbf{Theorem 4.2}. Assume the Conditions 3.1 and 4.1 are satisfied and 
\[
\frac{1}{q}+\frac{n}{r}=\frac{n}{2}-\gamma =\frac{1}{\tilde{q}}+\frac{n}{%
\tilde{r}}-2,\text{ }0\leq \alpha <1.
\]%
Let 
\[
\varphi \in \mathring{W}^{2,\gamma }\left( R^{n};H\left( A\right) \right) ,%
\text{ }\psi \in \mathring{W}^{2,\gamma -1}\left( R^{n};H\left( A\right)
\right) ,
\]%
\[
F\in L^{\tilde{q}^{\prime }}\left( \left[ 0,T\right] ;L^{\tilde{r}^{\prime
}}\left( R^{n};H\right) \right) 
\]
and let $u$ : $\left[ 0,T\right] \times R^{n}\rightarrow H\left( A\right) $
be a solution to $\left( 3.7-3.9\right) $. Then%
\[
\left\Vert A^{\alpha }u\right\Vert _{L^{q}\left( \left[ 0,T\right]
;L^{r}\left( R^{n};H\right) \right) }+\left\Vert A^{\alpha }u\right\Vert
_{C\left( \left[ 0,T\right] ;L^{2}\left( R^{n};H\right) \right) }+
\]%
\begin{equation}
\left\Vert A^{\alpha }\partial _{t}u\right\Vert _{C\left( \left[ 0,T\right] ;%
\mathring{W}^{2,\gamma -1}\left( R^{n};H\right) \right) }\lesssim \left\Vert
A\varphi \right\Vert _{\mathring{W}^{2,\gamma }\left( R^{n};H\right) }+ 
\tag{4.18}
\end{equation}%
\[
\left\Vert A\psi \right\Vert _{\mathring{W}^{2,\gamma -1}\left(
R^{n};H\right) }+\left\Vert F\right\Vert _{L^{\tilde{q}^{\prime }}\left( %
\left[ 0,T\right] ;L^{\tilde{r}^{\prime }}\left( R^{n};H\right) \right) }.
\]

\textbf{Proof. }By $\left( 3.18\right) -\left( 3.20\right) $ the solution of 
$\left( 3.7-3.9\right) $ can be expressed as .%
\begin{equation}
u\left( t,x\right) =S_{1}\left( t,x,A\right) \varphi \left( x\right)
+S_{2}\left( t,x,A\right) \psi \left( x\right) +\Phi \left( t,x,A\right) , 
\tag{4.19}
\end{equation}%
where 
\[
S_{1}\left( t,x,A\right) \varphi =F^{-1}D^{-1}\left( \xi \right) B_{1}\left(
t,\xi ,A\right) \hat{\varphi}\left( \xi \right) ,\text{ } 
\]

\begin{equation}
S_{2}\left( t,x,A\right) \psi =F^{-1}D^{-1}\left( \xi \right) B_{2}\left(
t,\xi ,A\right) \hat{\psi}\left( \xi \right) ,  \tag{4.20}
\end{equation}%
here

\[
B_{1}\left( t,\xi ,A\right) =\left[ C\left( t,\xi ,A\right) \left(
I-\dsum\limits_{k=1}^{m}\beta _{k}C\left( \lambda _{k},\xi ,A\right) \right)
\right. + 
\]%
\begin{equation}
\left. S\left( t,\xi ,A\right) \dsum\limits_{k=1}^{m}\beta _{k}A_{\xi
}S\left( \lambda _{k},\xi ,A\right) \right] ,  \tag{4.21}
\end{equation}%
\[
B_{2}\left( t,\xi ,A\right) =\left[ -C\left( t,\xi ,A\right)
\dsum\limits_{k=1}^{m}\alpha _{k}S\left( \lambda _{k},\xi ,A\right) \right.
+ 
\]%
\[
\left. S\left( t,\xi ,A\right) \left( I-\dsum\limits_{k=1}^{m}\alpha
_{k}C\left( \lambda _{k},\xi ,A\right) \right) \right] . 
\]

By the usual reduction using Littlewood-Paley theory we may assume that the
spatial Fourier transform of $\varphi $, $\psi ,$ $F$ and $u$ are all
localized in the annulus $\left\{ \left\vert \xi \right\vert \sim
2^{j}\right\} $ for some $j$ in a similar way as scular case (see Corollary
1.3 in $\left[ 22\right] $ and Lemma 5.1 of $[30]$ and the subsequent
discussion). The cases $r=\infty $ or $\tilde{r}=\infty $ can also be
treated by this argument, but the $H-$valued Lebesgue spaces $L_{x}^{r}$, $%
L_{x}^{\tilde{r}^{\prime }}$ must be replaced by their $H-$valued Besov
space counterparts. By the gap condition, the estimate is scale invariant,
and so we may assume $j=0.$ Now that frequency is localized, $\left( -\Delta
+A\right) ^{\frac{1}{2}}$ becomes an invertible smoothing operator, and we
may replace the Sobolev norms $\mathring{W}^{2,\gamma }\left( R^{n};H\right) 
$, $\mathring{W}^{2,\gamma -1}\left( R^{n};H\right) $ with the $L^{2}\left(
R^{n};H\right) $ norm. Combining these reductions with $\left( 4.19\right)
-\left( 4.21\right) $, we see that $\left( 4.18\right) $ will follow from
the estimates 
\[
\left\Vert S_{i\pm }\left( t,x,A\right) \varphi \right\Vert _{C\left(
L^{2}\left( R^{n};H\right) \right) }\lesssim \left\Vert \varphi \right\Vert
_{L^{2}\left( R^{n};H\right) }, 
\]%
\[
\left\Vert S_{i\pm }\left( t,x,A\right) \varphi \right\Vert
_{L_{t}^{q}\left( L_{x}^{r}\left( R^{n};H\right) \right) }\lesssim
\left\Vert \varphi \right\Vert _{L^{2}\left( R^{n};H\right) }, 
\]%
\[
\left\Vert S_{i\pm }\left( t,x,A\right) \psi \right\Vert _{C\left(
L^{2}\left( R^{n};H\right) \right) }\lesssim \left\Vert \psi \right\Vert
_{L^{2}\left( R^{n};H\right) }, 
\]%
\begin{equation}
\left\Vert S_{i\pm }\left( t,x,A\right) \psi \right\Vert _{L_{t}^{q}\left(
L_{x}^{r}\left( R^{n};H\right) \right) }\lesssim \left\Vert \psi \right\Vert
_{L^{2}\left( R^{n};H\right) },  \tag{4.22}
\end{equation}%
\[
\left\Vert \dint\limits_{t>s}S_{i\pm }\left( t,x,A\right) S_{i\pm }^{\ast
}\left( s,x,A\right) F\left( s\right) \right\Vert _{C\left( L^{2}\left(
R^{n};H\right) \right) }\lesssim \left\Vert F\right\Vert _{L_{t}^{q^{\prime
}}L_{x}^{\tilde{r}^{\prime }},} 
\]%
\[
\left\Vert \dint\limits_{t>s}S_{i\pm }\left( t,x,A\right) S_{i\pm }^{\ast
}\left( s,x,A\right) F\left( s\right) \right\Vert _{L_{t}^{q}L_{x}^{r}\left(
H\right) }\lesssim \left\Vert F\right\Vert _{L_{t}^{\tilde{q}^{\prime
}}L_{x}^{\tilde{r}^{\prime }}},\text{ }i=1,2, 
\]%
where the truncated wave evolution operators $S_{i\pm }\left( t,x,A\right) $
are given by%
\[
\hat{S}_{i\pm }\left( t,\xi ,A\right) f\left( \xi \right) =\chi _{\left[ 0,T%
\right] }\left( t\right) \beta \left( \xi \right) \hat{S}_{i}\left( t,\xi
,A\right) 
\]%
for some Littlewood-Paley cutoff function $\beta $\ supported on $\left\{
\left\vert \xi \right\vert \sim 1\right\} $. Apply Theorem 4.1 with all of
the above estimates $\left( 4.22\right) $ will follow from Theorem $4.1$
with $\Omega =R^{n}$, $Q=L^{2}(R^{n};H)$, once we show that operator
functions $S_{i\pm }\left( t,x,A\right) $ obey the energy estimate $(4.1)$
and the truncated decay estimate $(4.3)$. Consider first, the nonendpoint
case. By the method of $TT^{\ast }$ and in view of $\left( 4.20\right)
-\left( 4.21\right) $\ it will follow once we prove 
\begin{equation}
\left\Vert \dint\limits_{s<t}A^{\alpha }S_{i\pm }\left( t-s,x,A\right)
F\left( s\right) ds\right\Vert _{L_{t}^{q}L_{x}^{r}\left( H\right) }\lesssim
\left\Vert F\right\Vert _{L_{t}^{q^{\prime }}L_{x}^{r^{\prime }}\left(
H\right) }.  \tag{4.23}
\end{equation}

The energy estiamate $\left( 3.4\right) $:%
\[
\left\Vert U_{\Delta +A}\left( t\right) f\right\Vert _{L_{x}^{2}\left(
H\right) }\lesssim \left\Vert f\right\Vert _{L_{x}^{2}\left( H\right) } 
\]%
follows from Plancherel's theorem, the untruncated decay estimate%
\[
\left\Vert U_{\Delta }\left( t-s\right) f\right\Vert _{L_{x}^{\infty }\left(
H\right) }\lesssim \left\vert t-s\right\vert ^{-\frac{n}{2}}\left\Vert
f\right\Vert _{L_{x}^{1}\left( H\right) }, 
\]%
from the equality 
\[
U_{\Delta +A}\left( t\right) f=U_{\Delta }\left( t\right) U_{A}\left(
t\right) f, 
\]%
the explicit representation of the wave evolution operator%
\[
U_{\Delta }\left( t\right) f\left( x\right) =\left( 2\pi t\right) ^{-\frac{n%
}{2}}\dint\limits_{R^{n}}e^{\frac{\left\vert x-y\right\vert ^{2}}{2t}%
}f\left( y\right) dy 
\]%
and from the estimate $\left( 3.5\right) $. Due to properties of the
operator $A$, grops $U_{\Delta +A}\left( t\right) ,$ by $\left( 4.20\right)
-\left( 4.21\right) $ and by the dispersive estimate $(3.4)$ we have 
\[
\left\Vert A^{\alpha }\Phi _{i}\right\Vert _{E}\lesssim
\dint\limits_{s<t}\left\Vert A^{\alpha }S_{i\pm }\left( t-s,x,A\right)
ds\right\Vert _{B\left( H\right) }\left\Vert F\left( s\right) \right\Vert
_{H}ds\lesssim 
\]%
\[
\dint\limits_{\mathbb{R}}\left\vert t-s\right\vert ^{-n\left( \frac{1}{2}-%
\frac{1}{p}\right) -\alpha }\left\Vert F\left( s\right) \right\Vert _{H}ds, 
\]%
where 
\[
\Phi _{i}=\dint\limits_{s<t}S_{i\pm }\left( t-s,x,A\right) F\left( s\right)
ds\text{, }i=1,2. 
\]

Moreover, from above estimate by the Hardy-Littlewood-Sobolev inequality, we
get

\begin{equation}
\left\Vert A^{\alpha }\Phi _{i}\right\Vert _{L_{t}^{q}L_{x}^{r}\left(
R^{n+1};H\right) }\lesssim \left\Vert \dint\limits_{\mathbb{R}}\left\vert
t-s\right\vert ^{-n\left( \frac{1}{2}-\frac{1}{p}\right) -\alpha }\left\Vert
F\left( s\right) \right\Vert _{L_{x}^{r^{\prime }}\left( R^{n};H\right)
}ds\right\Vert _{L_{t}^{q}\left( \mathbb{R}\right) }\lesssim  \tag{4.24}
\end{equation}%
\[
\left\Vert F\right\Vert _{L_{t}^{q_{1}}L_{x}^{r^{\prime }}\left( H\right) }, 
\]%
where 
\[
\frac{1}{q_{1}}=\frac{1}{q}+\frac{1}{p}+\frac{1}{2}-\frac{\alpha }{n}. 
\]

The argument just presented also covers $\left( 4.24\right) $ in the case $q=%
\tilde{q},r=\tilde{r}$. It allows to consider the estimate in dualized form:%
\begin{equation}
\left\vert \dint \dint\limits_{s<t}\langle S_{i\pm }\left( t-s,x,A\right)
F\left( s\right) ,G\left( t\right) \rangle dsdt\right\vert \lesssim
\left\Vert F\right\Vert _{L_{t}^{q^{\prime }}L_{x}^{r^{\prime }}\left(
H\right) }\left\Vert G\right\Vert _{L_{t}^{\tilde{q}_{1}}L_{x}^{\tilde{r}%
^{\prime }}\left( H\right) }  \tag{4.25}
\end{equation}%
when 
\[
\frac{1}{\tilde{q}_{1}}=\frac{1}{\tilde{q}}+\frac{1}{\tilde{p}}+\frac{1}{2}-%
\frac{\alpha }{n}. 
\]%
The case $\tilde{q}=\infty ,$ $\tilde{r}=2$ follows from $\left( 4.25\right) 
$. Now, consider the endpoint case, i.e. $\left( q,r\right) =\left( 2,\frac{%
2n}{n-2}\right) $. It is suffices to show the following estimates%
\begin{equation}
\left\Vert A^{\alpha }S_{i\pm }\left( t-s,x,A\right) \varphi \right\Vert
_{L_{t}^{q}L_{x}^{r}\left( H\right) }\lesssim \left\Vert A\varphi
\right\Vert _{W^{s,2}\left( R^{n};H\right) },  \tag{4.26}
\end{equation}%
\[
\left\Vert A^{\alpha }S_{i\pm }\left( t-s,x,A\right) \varphi \right\Vert
_{C^{0}\left( L_{x}^{2}\left( H\right) \right) }\lesssim \left\Vert A\varphi
\right\Vert _{W^{s,2}\left( R^{n};H\right) }, 
\]%
\begin{equation}
\left\Vert A^{\alpha }S_{i\pm }\left( t-s,x,A\right) \psi \right\Vert
_{L_{t}^{q}L_{x}^{r}\left( H\right) }\lesssim \left\Vert A\psi \right\Vert
_{W^{s-1,2}\left( R^{n};H\right) },  \tag{4.27}
\end{equation}%
\[
\left\Vert A^{\alpha }S_{i\pm }\left( t-s,x,A\right) \psi \right\Vert
_{C^{0}\left( L_{x}^{2}\left( H\right) \right) }\lesssim \left\Vert A\psi
\right\Vert _{W^{s-1,2}\left( R^{n};H\right) }, 
\]

\begin{equation}
\left\Vert \dint\limits_{s<t}A^{\alpha }S_{i\pm }\left( t-s,x,A\right)
F\left( s\right) ds\right\Vert _{L_{t}^{q}L_{x}^{r}\left( H\right) }\lesssim
\left\Vert F\right\Vert _{L_{t}^{\tilde{q}^{\prime }}L_{x}^{\tilde{r}%
^{\prime }}\left( H\right) },  \tag{4.28}
\end{equation}%
\begin{equation}
\left\Vert \dint\limits_{s<t}A^{\alpha }S_{i\pm }\left( t-s,x,A\right)
F\left( s\right) ds\right\Vert _{C^{0}L_{x}^{2}\left( H\right) }\lesssim
\left\Vert F\right\Vert _{L_{t}^{q^{\prime }}L_{x}^{\tilde{r}^{\prime
}}\left( H\right) }.  \tag{4.29}
\end{equation}

Indeed, applying Theorem 4.1 with the energy estimate 
\[
\left\Vert S_{i\pm }\left( t-s,x,A\right) f\right\Vert _{L^{2}\left(
R^{n};H\right) }\lesssim \left\Vert f\right\Vert _{L^{2}\left(
R^{n};H\right) } 
\]%
which follows from Plancherel's theorem, the untruncated decay estimate $%
\left( 4.3\right) $ and by using of Lemma 4.1 we obtain the estimates $%
\left( 4.27\right) $ and $\left( 4.28\right) .$ Let us temporarily replace
the $C_{t}^{0}L_{x}^{2}\left( H\right) $ norm in estimates $\left(
4.26\right) $, $\left( 4.27\right) $ by the $L_{t}^{\infty }L_{x}^{2}\left(
H\right) .$ Then, all of the above the estimates will follow from Theorem
4.1, once we show that $S_{i\pm }\left( t,x,A\right) $\ obey the energy
estimate $\left( 4.1\right) $ and the truncated decay estimate $(4.2)$. The
estimate $\left( 4.1\right) $ is obtain immediate from Plancherel's theorem,
and $\left( 4.2\right) $ follows in a similar way as in $\left[ \text{31, p.
223-224}\right] $. To show that the quantity 
\[
G_{i}F\left( t\right) =\dint\limits_{s<t}A^{\alpha }S_{i\pm }\left(
t-s,x,A\right) F\left( s\right) ds\text{, }i=1,2 
\]%
is continuous in $L^{2}\left( R^{n};H\right) ,$ we use the the identity 
\[
G_{i}F\left( t+\varepsilon \right) =S_{i\pm }\left( \varepsilon ,x,A\right)
G_{i}F\left( t\right) +G_{i}\left( \chi _{\left[ t,t+\varepsilon \right]
}F\right) \left( t\right) , 
\]%
the continuity of $S_{i\pm }\left( \varepsilon ,x,A\right) $ as an operator
on $L^{2}\left( R^{n};H\right) $, and the fact that 
\[
\left\Vert \chi _{\left[ t,t+\varepsilon \right] }F\right\Vert _{L_{t}^{%
\tilde{q}^{\prime }}L_{x}^{\tilde{r}^{\prime }}\left( H\right) }\rightarrow 0%
\text{ as }\varepsilon \rightarrow 0. 
\]

From the estimates $\left( 4.26\right) -\left( 4.29\right) $ we obtain $%
\left( 4.18\right) $ for endpoint case.

\begin{center}
\textbf{5. Strichartz type estimates for solution} \textbf{to nonlinear wave
equation}
\end{center}

\bigskip\ For the Cauchy problem for scalar wave equation 
\begin{equation}
\partial _{t}^{2}u-\Delta u=F\left( u\right) ,\text{ }x\in R^{n},\text{ }%
t\in \left[ 0,T\right] ,  \tag{5.1}
\end{equation}%
\[
u\left( 0,x\right) =\varphi \left( .\right) \in \mathring{H}^{\gamma }\left(
R^{n}\right) ,\text{ }u_{t}\left( 0,x\right) =\psi \left( .\right) \in 
\mathring{H}^{\gamma -1}\left( R^{n}\right) , 
\]%
where nonlinearity $F\in C^{1}$ satisfies%
\[
F\left( u\right) =O\left( \left\vert u\right\vert ^{k}\right) ,\text{ }%
\left\vert u\right\vert \left\vert F_{u}\left( u\right) \right\vert \sim
\left\vert F\left( u\right) \right\vert \text{, }k>1. 
\]%
\ The question of how much regularity $\gamma =\gamma \left( k,n\right) $ is
needed to insure local well-posedness of this problem was addressed for
higher dimensions and nonlinearities in $[20]$; and then almost completely
answered in $[24]$.

Let%
\[
\text{ }X=L^{2}\left( R^{n};H\right) ,\text{ }Y=W^{2,2}\left( R^{n};H\left(
A\right) ,H\right) ,\text{ }H_{j}=\left( X,Y\right) _{\frac{1+2j}{4},2}\text{%
, }j=0,1. 
\]

\bigskip \textbf{Remark 5.1. }By using J. Lions-I. Petree result (see e.g $%
\left[ \text{38, \S\ 1.8.}\right] $) we obtain that the map $u\rightarrow
u^{\left( j\right) }\left( t_{0}\right) $, $t_{0}\in \left[ 0,T\right] $ is
continuous from $W^{2,2}\left( 0,T;X,Y\right) $ onto $H_{j}$ and there is a
constant $C_{1}$ such that 
\[
\left\Vert u^{\left( j\right) }\left( t_{0}\right) \right\Vert _{H_{j}}\leq
C_{1}\left\Vert u\right\Vert _{W^{2,2}\left( 0,T;X,Y\right) },\text{ }1\leq
p\leq \infty \text{.} 
\]

Consider the multipoint initial-value problem $\left( 1.1\right) -\left(
1.3\right) .$ By reasoning as $\left[ \text{22, Corollery 9.1}\right] $ we
prove the following result:

\textbf{Theorem 5.1. }Assume: (1) Conditions 3.1 and 4.1 are satisfied$;$

(2) the function $F:$ $H_{1}\rightarrow H$ is continuously differentiable
and obeys the power type estimates

\begin{equation}
F\left( u\right) =O\left( \left\Vert u\right\Vert ^{k}\right) ,\text{ }%
\left\Vert u\right\Vert \left\Vert F_{u}\left( u\right) \right\Vert \sim
\left\Vert F\left( u\right) \right\Vert \text{ }  \tag{5.2}
\end{equation}%
for some $k>1,$ where $F_{u}\left( u\right) $\ denotes the derivative of
operator function $F$ with respect to $u\in H$ and here $\left\Vert
u\right\Vert =\left\Vert u\right\Vert _{H};$

(3) $n\geq 4$, $\gamma =\frac{n-3}{2\left( n-1\right) }$, $k_{0}=\frac{%
\left( n+1\right) ^{2}}{\left( n-1\right) ^{2}+4};$

(4) $\varphi \in \mathring{W}^{2,\gamma }\left( R^{n};H\left( A\right)
\right) $, $\psi \in \mathring{W}^{2,\gamma -1}\left( R^{n};H\left( A\right)
\right) .$

Then for $k\geq k_{0}$ there is a $T>0$ depending only on 
\[
\left\Vert \varphi \right\Vert _{\mathring{W}^{2,\gamma }\left(
R^{n};H\left( A\right) \right) }+\left\Vert \psi \right\Vert _{\mathring{W}%
^{2,\gamma -1}\left( R^{n};H\left( A\right) \right) } 
\]%
and a unique weak solution $u$ to $\left( 1.1\right) -\left( 1.3\right) $
with 
\[
u\in L^{q_{0}}\left( \left[ 0,T\right] ;L^{r_{0}}\left( R^{n};H\left(
A\right) \right) \right) , 
\]%
where 
\[
q_{0}=\frac{2\left( n+1\right) }{n-3}\text{, }r_{0}=\frac{2\left(
n^{2}-1\right) }{\left( n^{2}-1\right) +4}\text{, }0\leq \alpha <1. 
\]

In addition, the solution satisfies 
\begin{equation}
u\in C\left( \left[ 0,T\right] ;\mathring{W}^{2,\gamma }\left( R^{n};H\left(
A\right) \right) \right) \cap C^{1}\left( \left[ 0,T\right] ;\mathring{W}%
^{2,\gamma -1}\left( R^{n};H\left( A\right) \right) \right)  \tag{5.3}
\end{equation}%
and depends continuously on the data.

\textbf{Proof. }We apply the standard fixed point argument in the space 
\[
V=V\left( T;M\right) =\left\{ u:u\in L^{q}\left( \left[ 0,T\right]
;L^{r}\left( R^{n};H\left( A\right) \right) \right) ,\right. \text{ } 
\]%
\[
\left. \left\Vert Au\right\Vert _{_{L_{t}^{q}L_{x}^{r}\left( H\right) }}\leq
M\right\} 
\]%
with $T$ and $M$ to be determined. Then we will used the estimate $\left(
4.18\right) $. By $(4.19)$, the problem of finding a solution $u$ of $%
(1.1)-\left( 1.3\right) $ is equivalent to finding a fixed point of the
mapping

\begin{equation}
S\left( u\right) =S_{1}\left( t,x,A\right) \varphi \left( x\right)
+S_{2}\left( t,x,A\right) \psi \left( x\right) +G\left( F\left( u\right)
\right) ,  \tag{5.4}
\end{equation}%
where $S_{i}\left( t,x,A\right) $, $i=1,2$ are operator function defined by $%
\left( 4.20\right) -\left( 4.21\right) $ and $G\left( F\left( u\right)
\right) $ defined by $\left( 3.20\right) $, where 
\begin{equation}
g_{1}\left( \xi \right) =\dsum\limits_{k=1}^{m}\alpha
_{k}\dint\limits_{0}^{\lambda _{k}}S\left( \lambda _{k}-\tau ,\xi ,A\right) 
\hat{F}\left( u\right) \left( \tau ,\xi \right) d\tau ,  \tag{5.5}
\end{equation}%
\[
g_{2}\left( \xi \right) =\dsum\limits_{k=1}^{m}\beta _{k}\left[ \frac{1}{2}%
\hat{F}\left( u\right) \left( \lambda _{k},\xi \right)
+\dint\limits_{0}^{\lambda _{k}}S\left( \lambda _{k}-\tau ,\xi ,A\right) 
\hat{F}\left( u\right) \left( \tau ,\xi \right) d\tau \right] . 
\]

\bigskip Accordingly, we will find $M,$ $T$ so that $S$ is a contraction on $%
V(T,M)$. It will suffice to show that for all $M$ there is a $T>0$ so that 
\begin{equation}
\left\Vert S\left( u\right) -S\left( \upsilon \right) \right\Vert _{V}\leq 
\frac{1}{2}\left\Vert u-\upsilon \right\Vert _{V}.  \tag{5.6}
\end{equation}

By $\left( 3.20\right) $ we get $G\left( F\left( 0\right) \right) =0.$ It
implies 
\[
S\left( 0\right) =S_{1}\left( t,x,A\right) \varphi \left( x\right)
+S_{2}\left( t,x,A\right) \psi \left( x\right) . 
\]
So $S\left( 0\right) $ is finite by $\left( 4.18\right) $ applied to the
homogeneous problem and by using the properties of functions $S_{1}$, $S_{2}$
with relation of the operator $A,$ i.e., the fact that $S$: $V\rightarrow V$
follows by picking $M$ large enough so 
\begin{equation}
\left\Vert S\left( 0\right) \right\Vert _{V}\leq \frac{M}{2}.  \tag{5.7}
\end{equation}

Again in view of $\left( 4.18\right) $ we have 
\begin{equation}
\left\Vert S\left( u\right) -S\left( \upsilon \right) \right\Vert
_{V}=G\left( F\left( u\right) -F\left( \upsilon \right) \right) \lesssim
\left\Vert F\left( u\right) -F\left( \upsilon \right) \right\Vert _{L_{t}^{%
\tilde{q}^{\prime }}L_{x}^{\tilde{r}^{\prime }}\left( H\right) }.  \tag{5.8}
\end{equation}

The assumptions $(5.2)$ give%
\[
\left\Vert F\left( u\right) -F\left( \upsilon \right) \right\Vert
_{H}=\left\Vert \dint\limits_{0}^{1}\frac{d}{d\lambda }F\left( \lambda
u+\left( 1-\lambda \right) \upsilon \right) d\lambda \right\Vert _{H}= 
\]%
\[
\left\Vert \dint\limits_{0}^{1}\left( u-\upsilon \right) \nabla F\left(
\lambda u+\left( 1-\lambda \right) \upsilon \right) d\lambda \right\Vert
_{H}\lesssim \left\Vert u-\upsilon \right\Vert _{H}\left( \left\Vert
u\right\Vert _{H}+\left\Vert \upsilon \right\Vert _{H}\right) ^{k-1}. 
\]

Using this in $\left( 5.8\right) $ gives 
\begin{equation}
\left\Vert S\left( u\right) -S\left( \upsilon \right) \right\Vert
_{V}\lesssim \left\Vert \left\Vert u-\upsilon \right\Vert _{H}\left(
\left\Vert u\right\Vert _{H}+\left\Vert \upsilon \right\Vert _{H}\right)
^{k-1}\right\Vert _{L_{t}^{\tilde{q}^{\prime }}L_{x}^{\tilde{r}^{\prime }}}.
\tag{5.9}
\end{equation}

Moreover, by the generalized H\"{o}lder inequality we have 
\begin{equation}
\left\Vert \left\Vert u-\upsilon \right\Vert _{H}\left( \left\Vert
u\right\Vert _{H}+\left\Vert \upsilon \right\Vert _{H}\right)
^{k-1}\right\Vert _{L_{t}^{\tilde{q}^{\prime }}L_{x}^{\tilde{r}^{\prime
}}}\leq  \tag{5.10}
\end{equation}%
\[
\left\Vert u-\upsilon \right\Vert _{L_{t}^{q}L_{x}^{r}\left( H\right)
}\left\Vert \left( \left\Vert u\right\Vert _{H}+\left\Vert \upsilon
\right\Vert _{H}\right) ^{k-1}\right\Vert
_{L_{t}^{q/k-1}L_{x}^{r/k-1}}\left\Vert \chi _{\left[ 0,T\right]
}\right\Vert _{L_{t}^{p}L_{x}^{\infty }\left( H\right) }, 
\]

where $1\leq p<\infty $ is chosen so that 
\[
\frac{1}{\tilde{q}^{\prime }}=\frac{1}{q}+\frac{1}{q/\left( k-1\right) }+%
\frac{1}{p}\text{, }\frac{1}{\tilde{r}^{\prime }}=\frac{1}{r}+\frac{1}{%
r/\left( k-1\right) }. 
\]

By the assumptions on $u$, $\upsilon $ the estimate $(5.10)$ simplifies to 
\begin{equation}
\left\Vert \left( \left\Vert u\right\Vert _{H}+\left\Vert \upsilon
\right\Vert _{H}\right) ^{k-1}\right\Vert _{L_{t}^{\tilde{q}^{\prime
}}L_{x}^{\tilde{r}^{\prime }}}\lesssim T^{\frac{1}{p}}M^{k-1}\left\Vert
u-\upsilon \right\Vert _{V}.  \tag{5.11}
\end{equation}

Thus if we choose $T$ so that $T^{\frac{1}{p}}M^{k-1}<1$, then $(5.9)$ and $%
(5.11)$ give the desired contraction $\left( 5.6\right) .$

To obtain the regularity $(5.3)$ for $u$ we apply $(5.11)$ with $\upsilon =0$
to obtain%
\[
\left\Vert S\left( u\right) \right\Vert _{L_{t}^{\tilde{q}^{\prime }}L_{x}^{%
\tilde{r}^{\prime }}\left( H\right) }\leq T^{\frac{1}{p}}M^{k-1}\left\Vert
u\right\Vert _{V}<\infty , 
\]%
and $(5.3)$ follows from $\left( 4.18\right) $.

Finally, we need to show uniqueness. Suppose that we have two solutions $u,$ 
$\upsilon $ to $\left( 1.1\right) -\left( 1.3\right) $ for time $[0,T^{\ast
}]$ such that 
\[
\left\Vert Au\right\Vert _{_{L^{q}\left( \left[ 0,T^{\ast }\right]
;L^{r}\left( R^{n};H\right) \right) }}\leq M,\text{ }\left\Vert A\upsilon
\right\Vert _{_{L^{q}\left( \left[ 0,T^{\ast }\right] ;L^{r}\left(
R^{n};H\right) \right) }}\leq M
\]%
for some $M.$ Choose $0<T\leq T^{\ast }$ such that $T^{\frac{1}{p}}M^{k-1}<1.
$ By the above arguments $(5.6)$ holds, which implies that $u=\upsilon $ for
time $[0,T]$. Since $T$ depends only on $M$, we may iterate this argument
and obtain $u=\upsilon $ for all times $[0,T^{\ast }]$.

\begin{center}
\textbf{6.The exsistence and uniquness for the system of wave equation }
\end{center}

Consider at first, the multipoint Cauchy problem for linear system of wave
equations 
\begin{equation}
\partial _{t}^{2}u_{m}-\Delta
u_{m}+\sum\limits_{j=1}^{N}a_{mj}u_{j}=F_{j}\left( t,x\right) ,\text{ }t\in %
\left[ 0,T\right] \text{, }x\in R^{n},  \tag{6.1}
\end{equation}%
\[
u_{m}\left( 0,x\right) =\varphi _{m}\left( x\right)
+\dsum\limits_{k=1}^{m}\alpha _{k}u_{m}\left( \lambda _{k},x\right) ,\text{
for a.e. }x\in R^{n}, 
\]%
\[
\partial _{t}u_{m}\left( 0,x\right) =\psi _{m}\left( x\right)
+\dsum\limits_{k=1}^{m}\beta _{k}\partial _{t}u_{m}\left( \lambda
_{k},x\right) ,\text{ for a.e. }x\in R^{n}, 
\]%
where $u=\left( u_{1},u_{2},...,u_{N}\right) ,$ $u_{j}=u_{j}\left(
t,x\right) ,$ $a_{mj}$ are complex numbers. Let $l_{2}=l_{2}\left( N\right) $
and $l_{2}^{s}=l_{2}^{s}\left( N\right) ,$ $N\in \mathbb{N}$ (see $\left[ 
\text{38, \S\ 1.18}\right] $). Let $A$ be the operator in $l_{2}\left(
N\right) $ defined by%
\[
\text{ }D\left( A\right) =\left\{ u=\left\{ u_{j}\right\} \in l_{2}\left(
N\right) ,\text{ }\left\Vert Au\right\Vert _{l_{2}\left( N\right) }=\left(
\sum\limits_{m,j=1}^{N}\left\vert a_{mj}u_{j}\right\vert ^{2}\right) ^{\frac{%
1}{2}}<\infty \right\} , 
\]

\[
A=\left[ a_{mj}\right] \text{, }a_{mj}=a_{jm},\text{ }s>0,\text{ }%
m,j=1,2,...,N,\text{ }N\in \mathbb{N}. 
\]

It is clear that for $N=m<\infty $ the space $l_{2}\left( N\right) $ conside
with the finite dimensional vector space $\mathbb{C}^{m}.$

\ From Theorem 4.2 we obtain the following result

\bigskip \textbf{Theorem 6.1. }Assume the Conditions 4.1 are hold. Let 
\[
\left\vert \alpha _{k}+\beta _{k}\right\vert
>0,\dsum\limits_{k,j=1}^{m}\alpha _{k}\beta _{j}\neq 0,\text{ }0\leq \alpha
<1
\]%
and 
\[
\frac{1}{q}+\frac{n}{r}=\frac{n}{2}-\gamma =\frac{1}{\tilde{q}}+\frac{n}{%
\tilde{r}}-2.
\]
Suppose $n\geq 1$ and 
\[
\varphi \in \mathring{W}^{\gamma ,2}\left( R^{n};D\left( A\right) \right) ,%
\text{ }\psi \in \mathring{W}^{\gamma -1,2}\left( R^{n};D\left( A\right)
\right) ,\text{ }
\]%
\[
F\in L^{\tilde{q}^{\prime }}\left( \left[ 0,T\right] ;L^{\tilde{r}^{\prime
}}\left( R^{n};l_{2}\left( N\right) \right) \right) .
\]%
Let $u$ : $\left[ 0,T\right] \times R^{n}\rightarrow l_{2}\left( N\right) $
be a solution to $\left( 6.1\right) $. Then%
\begin{equation}
\left\Vert A^{\alpha }u\right\Vert _{L^{q}\left( \left[ 0,T\right]
;L^{r}\left( R^{n};l_{2}\left( N\right) \right) \right) }+\left\Vert
A^{\alpha }u\right\Vert _{C\left( \left[ 0,T\right] ;L^{2}\left(
R^{n};l_{2}\left( N\right) \right) \right) }+  \tag{6.2}
\end{equation}%
\[
\left\Vert A^{\alpha }\partial _{t}u\right\Vert _{C\left( \left[ 0,T\right] ;%
\mathring{W}^{2,\gamma -1}\left( R^{n};l_{2}\left( N\right) \right) \right)
}\lesssim \left\Vert A\varphi \right\Vert _{\mathring{W}^{2,\gamma }\left(
R^{n};l_{2}\left( N\right) \right) }+
\]%
\[
\left\Vert A\psi \right\Vert _{\mathring{W}^{2,\gamma -1}\left(
R^{n};l_{2}\left( N\right) \right) }+\left\Vert F\right\Vert _{L^{\tilde{q}%
^{\prime }}\left( \left[ 0,T\right] ;L^{\tilde{r}^{\prime }}\left(
R^{n};l_{2}\left( N\right) \right) \right) }.
\]

\ \textbf{Proof.} It is easy to see that $A$ is a symmetric operator in $%
l_{2}\left( N\right) $ and other conditions of Theorem 4.2 are satisfied.
Hence, from Theorem 4.2 we obtain the conculision.

Consider now, the Cauchy problem $\left( 1.10\right) $. Let $A$ be the
operator in $l_{2}\left( N\right) $ defined by

\[
\text{ }D\left( A\right) =\left\{ u=\left\{ u_{j}\right\} \in l_{2}^{s},%
\text{ }s>0\text{ }\left\Vert u\right\Vert _{l_{2}^{s}\left( N\right)
}=\left( \sum\limits_{j=1}^{N}\left\vert 2^{sj}u_{j}\right\vert ^{2}\right)
^{\frac{1}{2}}<\infty \right\} , 
\]%
\[
A=\left[ a_{mj}\right] \text{, }a_{mj}=a_{jm},\text{ }s>0,\text{ }%
m,j=1,2,...,N,\text{ }N\in \mathbb{N}, 
\]

\ \ Let 
\[
\text{ }X\left( N\right) =L^{2}\left( R^{n};l_{2}\left( N\right) \right) ,%
\text{ }Y\left( N\right) =W^{2,2}\left( R^{n};l_{2}^{s}\left( N\right)
,l_{2}\left( N\right) \right) ,\text{ } 
\]%
\[
H_{j}\left( N\right) =\left( X\left( N\right) ,Y\left( N\right) \right) _{%
\frac{1+2j}{4},2}\text{, }j=0,1, 
\]%
where $H_{j}\left( N\right) $ denote the real interpolation spaces between $%
X\left( N\right) $ and $Y\left( N\right) .$

\textbf{Remark 6.1. } It is known that (see e.g. $\left[ \text{38, \S\ 1.18}%
\right] $ the real interpolation spaces $\left( l_{2}^{s}\left( N\right)
,l_{2}\left( N\right) \right) _{\theta ,2}$, $\theta \in \left( 0,1\right) $
between $l_{2}^{s}\left( N\right) $ and $l_{2}\left( N\right) $ defined as 
\[
\left( l_{2}^{s}\left( N\right) ,l_{2}\left( N\right) \right) _{\theta
,2}=l_{2}^{s\left( 1-\theta \right) }\left( N\right) . 
\]
So, it can be shown that\textbf{\ }%
\[
H_{j}\left( N\right) =W^{2\left( 1-\theta \right) ,2}\left(
R^{n};l_{2}^{s\left( 1-\theta \right) }\left( N\right) ,l_{2}\left( N\right)
\right) . 
\]

We obtain from Theorem 5.1 the following result

\textbf{Theorem 6.2. }Assume: (1) the function $F:$ $l_{2}^{s/4}\left(
N\right) \rightarrow l_{2}\left( N\right) $ is continuously differentiable
and obeys the power type estimates

\[
F\left( u\right) =O\left( \left\Vert u\right\Vert _{l_{2}}^{k}\right) ,\text{
}\left\Vert u\right\Vert _{l_{2}}\left\Vert F_{u}\left( u\right) \right\Vert
_{l_{2}}\sim \left\Vert F\left( u\right) \right\Vert _{l_{2}}\text{ } 
\]%
for some $k>1;$

(2) $n\geq 4$, $\gamma =\frac{n-3}{2\left( n-1\right) }$, $k=\frac{\left(
n+1\right) ^{2}}{\left( n-1\right) ^{2}+4};$

(3) $\varphi \in \mathring{W}^{2,\gamma }\left( R^{n};l_{2}^{s}\left(
N\right) \right) $, $\psi \in \mathring{W}^{2,\gamma -1}\left(
R^{n};l_{2}^{s}\left( N\right) \right) ;$

(4) assume the Conditions 4.1 are hold and 
\[
\left\vert \alpha _{k}+\beta _{k}\right\vert
>0,\dsum\limits_{k,j=1}^{m}\alpha _{k}\beta _{j}\neq 0, 
\]%
\[
\frac{1}{q}+\frac{n}{r}=\frac{n}{2}-\gamma =\frac{1}{\tilde{q}}+\frac{n}{%
\tilde{r}}-2. 
\]

Then there is a $T>0$ depending only on 
\[
\left\Vert \varphi \right\Vert _{\mathring{W}^{2,\gamma }\left(
R^{n};l_{2}^{s}\left( N\right) \right) }+\left\Vert \psi \right\Vert _{%
\mathring{W}^{2,\gamma -1}\left( R^{n};l_{2}^{s}\left( N\right) \right) }
\]%
and a unique weak solution $u$ to $\left( 1.12\right) -\left( 1.9\right) $
with 
\[
u\in L^{q_{0}}\left( \left[ 0,T\right] ;L^{r_{0}}\left(
R^{n};l_{2}^{s}\left( N\right) \right) \right) ,
\]%
where 
\[
q_{0}=\frac{2\left( n+1\right) }{n-3}\text{, }r_{0}=\frac{2\left(
n^{2}-1\right) }{\left( n^{2}-1\right) +4}.
\]

In addition, the solution satisfies 
\[
u\in C\left( \left[ 0,T\right] ;\mathring{W}^{2,\gamma }\left(
R^{n};l_{2}^{s}\left( N\right) \right) \right) \cap C^{1}\left( \left[ 0,T%
\right] ;\mathring{W}^{2,\gamma -1}\left( R^{n};l_{2}^{s}\left( N\right)
\right) \right) 
\]%
and depends continuously on the data.

\ \textbf{Proof.} It is easy to see that $A$ is a symmetric operator in $%
l_{2}$ and other conditions of Theorem 5.1 are satisfied. Hence, from Teorem
5.1 we obtain the conculision.

\begin{center}
\textbf{7.The exsistence and uniquness of solution to anisotropic wave
equation}\ \ \ \ \ \ \ \ \ \ \ \ \ \ \ \ \ \ \ \ \ \ \ \ \ \ \ \ \ \ \ \ \ \
\ \ \ \ \ \ \ \ \ \ \ \ \ \ 
\end{center}

Let $\Omega =R^{n}\times G$, $G\subset R^{d},$ $d\geq 2$ is a bounded domain
with $\left( d-1\right) $-dimensional boundary $\partial G$. Consider at
first, the multipoint mixed problem for the following wave equation

\begin{equation}
\partial _{t}^{2}u-\Delta _{x}u+\sum\limits_{\left\vert \alpha \right\vert
\leq 2l}a_{\alpha }\left( y\right) D_{y}^{\alpha }u=F\left( t,x\right) ,%
\text{ }  \tag{7.1}
\end{equation}%
\[
\text{ }x\in R^{n},\text{ }y\in G,\text{ }t\in \left[ 0,T\right] ,\text{ }%
p\geq 0, 
\]

\begin{equation}
B_{j}u=\sum\limits_{\left\vert \beta \right\vert \leq l_{j}}\ b_{j\beta
}\left( y\right) D_{y}^{\beta }u=0\text{, }x\in R^{n},\text{ }y\in \partial
G,\text{ }j=1,2,...,m,  \tag{7.2}
\end{equation}%
\begin{equation}
u\left( 0,x,y\right) =\varphi \left( x,y\right)
+\dsum\limits_{k=1}^{m}\alpha _{k}u\left( \lambda _{k},x,y\right) ,\text{
for a.e. }x\in R^{n},\text{ }y\in G,  \tag{7.3}
\end{equation}%
\[
\partial _{t}u\left( 0,x,y\right) =\psi \left( x,y\right)
+\dsum\limits_{k=1}^{m}\beta _{k}\partial _{t}u\left( \lambda
_{k},x,y\right) ,\text{ for a.e. }x\in R^{n},\text{ }y\in G, 
\]%
where $u=u\left( t,x,y\right) $ is a solution, $a_{\alpha },$ $b_{j\beta }$
are the complex valued functions, $\lambda =\pm 1,$ $\alpha =\left( \alpha
_{1},\alpha _{2},...,\alpha _{d}\right) $, $\beta =\left( \beta _{1},\beta
_{2},...,\beta _{d}\right) ,$ $\mu _{i}<2l$ and 
\[
D_{x}^{k}=\frac{\partial ^{k}}{\partial x^{k}},\text{ }D_{j}=-i\frac{%
\partial }{\partial y_{j}},\text{ }D_{y}=\left( D_{1,}...,D_{d}\right) ,%
\text{ }y=\left( y_{1},...,y_{d}\right) . 
\]

$\ $

\bigskip Let%
\[
\xi ^{\prime }=\left( \xi _{1},\xi _{2},...,\xi _{d-1}\right) \in R^{d-1},%
\text{ }\alpha ^{\prime }=\left( \alpha _{1},\alpha _{2},...,\alpha
_{d-1}\right) \in Z^{d-1},\text{ } 
\]%
\[
\text{ }A\left( y_{0},\xi ^{\prime },D_{y}\right) =\sum\limits_{\left\vert
\alpha ^{\prime }\right\vert +j\leq 2l}a_{\alpha ^{\prime }}\left(
y_{0}\right) \xi _{1}^{\alpha _{1}}\xi _{2}^{\alpha _{2}}...\xi
_{d-1}^{\alpha _{d-1}}D_{y}^{j}\text{ for }y_{0}\in \bar{G} 
\]%
\[
B_{j}\left( y_{0},\xi ^{\prime },D_{y}\right) =\sum\limits_{\left\vert \beta
^{\prime }\right\vert +j\leq l_{j}}b_{j\beta ^{\prime }}\left( y_{0}\right)
\xi _{1}^{\beta _{1}}\xi _{2}^{\beta _{2}}...\xi _{d-1}^{\beta
_{d-1}}D_{y}^{j}\text{ for }y_{0}\in \partial G. 
\]

For $\Omega =R^{n}\times G,$ $s\in \mathbb{R}$ and $l\in \mathbb{N}$ let $%
\mathring{W}^{s,l,p}\left( \Omega \right) =\mathring{W}^{s,l,p}\left( \Omega
;\mathbb{C}\right) .$

From Theorem 4.2 we obtain the following result

\bigskip \textbf{Theorem 7.1. }Assume the following conditions be satisfied:

\bigskip (1) $G\in C^{2}$, $a_{\alpha }\in C\left( \bar{G}\right) $ for each 
$\left\vert \alpha \right\vert =2l$ and $a_{\alpha }\in L_{\infty }\left(
G\right) $ for each $\left\vert \alpha \right\vert <2l$;

(2) $b_{j\beta }\in C^{2l-l_{j}}\left( \partial G\right) $ for each $j$, $%
\beta $ and $\ l_{j}<2l$, $\sum\limits_{j=1}^{l}b_{j\beta }\left( y^{\prime
}\right) \sigma _{j}\neq 0,$ for $\left\vert \beta \right\vert =m_{j},$ $%
y^{^{\shortmid }}\in \partial G,$ where $\sigma =\left( \sigma _{1},\sigma
_{2},...,\sigma _{d}\right) \in R^{d}$ is a normal to $\partial G$ $;$

(3) for $y\in \bar{G}$, $\xi \in R^{d}$, $\mu \in S\left( \varphi
_{0}\right) $ for $0\leq \varphi _{0}<\pi $, $\left\vert \xi \right\vert
+\left\vert \mu \right\vert \neq 0$ let $\mu +$ $\sum\limits_{\left\vert
\alpha \right\vert =2l}a_{\alpha }\left( y\right) \xi ^{\alpha }\neq 0$;

(4) for each $y_{0}\in \partial G$ local BVP in local coordinates
corresponding to $y_{0}$:%
\[
\mu +A\left( y_{0},\xi ^{\prime },D_{y}\right) \vartheta \left( y\right) =0, 
\]

\[
B_{j}\left( y_{0},\xi ^{\prime },D_{y}\right) \vartheta \left( 0\right)
=h_{j}\text{, }j=1,2,...,l 
\]%
has a unique solution $\vartheta \in C_{0}\left( \mathbb{R}_{+}\right) $ for
all $h=\left( h_{1},h_{2},...,h_{d}\right) \in \mathbb{C}^{d}$ and for $\xi
^{\prime }\in R^{d-1};$

(5) assume the Conditions 4.1 are hold and 
\[
\left\vert \alpha _{k}+\beta _{k}\right\vert
>0,\dsum\limits_{k,j=1}^{m}\alpha _{k}\beta _{j}\neq 0,\text{ }0\leq \alpha
<1,
\]%
\[
\frac{1}{q}+\frac{n}{r}=\frac{n}{2}-\gamma =\frac{1}{\tilde{q}}+\frac{n}{%
\tilde{r}}-2;
\]

(6) $n\geq 1$ and 
\[
\varphi \in \mathring{W}^{\gamma ,2}\left( R^{n};W^{2l,2}\left( G\right)
\right) ,\text{ }\psi \in \mathring{W}^{\gamma -1,2}\left(
R^{n};W^{2l,2}\left( G\right) \right) ,\text{ }
\]%
\[
F\in L^{\tilde{q}^{\prime }}\left( \left[ 0,T\right] ;L^{\tilde{r}^{\prime
}}\left( R^{n};L^{2}\left( G\right) \right) \right) .
\]
Let $u$ : $\left[ 0,T\right] \times R^{n}\rightarrow L^{2}\left( G\right) $
be a solution to $\left( 5.1\right) $. Then%
\[
\left\Vert A^{\alpha }u\right\Vert _{L^{q}\left( \left[ 0,T\right]
;L^{r}R^{n};L^{2}\left( G\right) \right) }+\left\Vert A^{\alpha
}u\right\Vert _{C\left( \left[ 0,T\right] ;L^{2}\left( R^{n};L^{2}\left(
G\right) \right) \right) }+
\]%
\[
\left\Vert A^{\alpha }\partial _{t}u\right\Vert _{C\left( \left[ 0,T\right] ;%
\mathring{W}^{2,\gamma -1}R^{n};L^{2}\left( G\right) \right) }\lesssim
\left\Vert A\varphi \right\Vert _{\mathring{W}^{2,\gamma }\left(
R^{n};L^{2}\left( G\right) \right) }+
\]%
\[
\left\Vert A\psi \right\Vert _{\mathring{W}^{2,\gamma -1}\left(
R^{n};L^{2}\left( G\right) \right) }+\left\Vert F\right\Vert _{L^{\tilde{q}%
^{\prime }}\left( \left[ 0,T\right] ;L^{\tilde{r}^{\prime
}}R^{n};L^{2}\left( G\right) \right) }.
\]

\textbf{Proof. }Let us consider the operator $A$ in $H=L^{2}\left( G\right) $
that are defined by 
\[
D\left( A\right) =\left\{ u\in W^{2l,2}\left( G\right) \text{, }B_{j}u=0,%
\text{ }j=1,2,...,l\text{ }\right\} ,\ Au=\sum\limits_{\left\vert \alpha
\right\vert \leq 2l}a_{\alpha }\left( y\right) D_{y}^{\alpha }u\left(
y\right) . 
\]

Then the problem $\left( 7.1\right) -\left( 7.3\right) $ can be rewritten as
the problem $\left( 3.7\right) -\left( 3.9\right) $, where $u\left( x\right)
=u\left( x,.\right) ,$ $f\left( x\right) =f\left( x,.\right) $,\ $x\in R^{n}$
are the functions with values in\ $H=L^{2}\left( G\right) $. By virtue of $%
\left[ \text{8, Theorem 8.2}\right] ,$ operator $A$ is absolute positive in $%
L^{2}\left( G\right) $. Moreover, in view of (1)-(6) all conditons of
Theorem 4.2 are hold. Then Theorem 4.2 implies the assertion.

Consider now, the multipoint mixed problem for nonlinear wave equation

\begin{equation}
\partial _{t}^{2}u+\Delta _{x}u+\sum\limits_{\left\vert \alpha \right\vert
\leq 2l}a_{\alpha }\left( y\right) D_{y}^{\alpha }u=F\left( u\right) ,\text{ 
}  \tag{7.4}
\end{equation}%
\[
\text{ }x\in R^{n},\text{ }y\in G,\text{ }t\in \left[ 0,T\right] ,\text{ }%
p\geq 0,
\]

\begin{equation}
B_{j}u=\sum\limits_{\left\vert \beta \right\vert \leq l_{j}}\ b_{j\beta
}\left( y\right) D_{y}^{\beta }u=0\text{, }x\in R^{n},\text{ }y\in \partial
G,\text{ }j=1,2,...,l,  \tag{7.5}
\end{equation}%
\begin{equation}
u\left( 0,x,y\right) =\varphi \left( x,y\right)
+\dsum\limits_{k=1}^{m}\alpha _{k}u\left( \lambda _{k},x,y\right) ,\text{
for a.e. }x\in R^{n},\text{ }y\in G  \tag{7.6}
\end{equation}%
\[
\partial _{t}u\left( 0,x,y\right) =\psi \left( x,y\right)
+\dsum\limits_{k=1}^{m}\beta _{k}\partial _{t}u\left( \lambda
_{k},x,y\right) ,\text{ for a.e. }x\in R^{n}\text{, }y\in G. 
\]

From Theorem 5.1 we obtain

\textbf{Theorem 7.2}. Assume the following conditions be satisfied:

\bigskip (1) $G\in C^{2}$, $a_{\alpha }\in C\left( \bar{G}\right) $ for each 
$\left\vert \alpha \right\vert =2l$ and $a_{\alpha }\in L_{\infty }\left(
G\right) $ for each $\left\vert \alpha \right\vert <2m$;

(2) $b_{j\beta }\in C^{2m-m_{j}}\left( \partial G\right) $ for each $j$, $%
\beta $ and $\ m_{j}<2l$, $\sum\limits_{j=1}^{l}b_{j\beta }\left( y^{\prime
}\right) \sigma _{j}\neq 0,$ for $\left\vert \beta \right\vert =l_{j},$ $%
y^{^{\shortmid }}\in \partial G,$ where $\sigma =\left( \sigma _{1},\sigma
_{2},...,\sigma _{d}\right) \in R^{d}$ is a normal to $\partial G$ $;$

(3) for $y\in \bar{G}$, $\xi \in R^{d}$, $\mu \in S\left( \varphi
_{0}\right) $ for $0\leq \varphi _{0}<\pi $, $\left\vert \xi \right\vert
+\left\vert \mu \right\vert \neq 0$ let $\mu +$ $\sum\limits_{\left\vert
\alpha \right\vert =2d}a_{\alpha }\left( y\right) \xi ^{\alpha }\neq 0$;

(4) for each $y_{0}\in \partial G$ local BVP in local coordinates
corresponding to $y_{0}$:%
\[
\mu +A\left( y_{0},\xi ^{\prime },D_{y}\right) \vartheta \left( y\right) =0, 
\]

\[
B_{j}\left( y_{0},\xi ^{\prime },D_{y}\right) \vartheta \left( 0\right)
=h_{j}\text{, }j=1,2,...,l 
\]%
has a unique solution $\vartheta \in C_{0}\left( \mathbb{R}_{+}\right) $ for
all $h=\left( h_{1},h_{2},...,h_{d}\right) \in \mathbb{C}^{d}$ and for $\xi
^{\prime }\in R^{d-1};$

(5) assume the Condition 4.1 are hold and and 
\[
\left\vert \alpha _{k}+\beta _{k}\right\vert
>0,\dsum\limits_{k,j=1}^{m}\alpha _{k}\beta _{j}\neq 0,
\]%
\[
\frac{1}{q}+\frac{n}{r}=\frac{n}{2}-\gamma =\frac{1}{\tilde{q}}+\frac{n}{%
\tilde{r}}-2;
\]

(6) the function $F:$ $W^{\frac{l}{2},2}\left( G\right) \rightarrow
L^{2}\left( G\right) $ is continuously differentiable and obeys the power
type estimates

\[
F\left( u\right) =O\left( \left\Vert u\right\Vert _{L^{2}\left( G\right)
}^{k}\right) ,\text{ }\left\Vert u\right\Vert _{L^{2}\left( G\right)
}\left\Vert F_{u}\left( u\right) \right\Vert _{L^{2}\left( G\right) }\sim
\left\Vert F\left( u\right) \right\Vert _{L^{2}\left( G\right) }\text{ }
\]%
for some $k>1;$

\bigskip (7) $n\geq 4$, $\gamma =\frac{n-3}{2\left( n-1\right) }$, $k=\frac{%
\left( n+1\right) ^{2}}{\left( n-1\right) ^{2}+4};$

(8) $\varphi \in \mathring{W}^{2,\gamma }\left( R^{n};W^{2l,2}\left(
G\right) \right) $ and $\psi \in \mathring{W}^{2,\gamma -1}\left(
R^{n};W^{2l,2}\left( G\right) \right) .$

\bigskip Then there is a $T>0$ depending only on 
\[
\left\Vert \varphi \right\Vert _{\mathring{W}^{2,\gamma }\left(
R^{n};W^{2l,2}\left( G\right) \right) }+\left\Vert \psi \right\Vert _{%
\mathring{W}^{2,\gamma -1}\left( R^{n};W^{2l,2}\left( G\right) \right) }
\]%
and a unique weak solution $u$ to $\left( 7.4.1\right) -\left( 7.6\right) $
with 
\[
u\in L^{q_{0}}\left( \left[ 0,T\right] ;L^{r_{0}}\left( R^{n};W^{2l,2}\left(
G\right) \right) \right) ,
\]%
where 
\[
q_{0}=\frac{2\left( n+1\right) }{n-3}\text{, }r_{0}=\frac{2\left(
n^{2}-1\right) }{\left( n^{2}-1\right) +4}.
\]

In addition, the solution satisfies 
\[
u\in C\left( \left[ 0,T\right] ;\mathring{W}^{2,\gamma }\left(
R^{n};W^{2l,2}\left( G\right) \right) \right) \cap C^{1}\left( \left[ 0,T%
\right] ;\mathring{W}^{2,\gamma -1}\left( R^{n};W^{2l,2}\left( G\right)
\right) \right) 
\]%
and depends continuously on the data.

\ \textbf{Proof. }The problem $\left( 7.4\right) -\left( 7.6\right) $ can be
rewritten as the problem $\left( 1.1\right) $, where $u\left( x\right)
=u\left( x,.\right) ,$ $f\left( x\right) =f\left( x,.\right) $,\ $x\in R^{n}$
are the functions with values in\ $H=L^{2}\left( G\right) $. By virtue of $%
\left[ \text{8, Theorem 8.2}\right] ,$ operator $A+\mu $ is absolute
positive in $L^{2}\left( G\right) $ for sufficiently large $\mu >0$.
Moreover, in view of (1)-(8) all conditons of Theorem 5.1 are hold. Then
Theorem 5.1 implies the assertion.

\begin{center}
\textbf{8.} \textbf{The Wentzell-Robin type mixed problem for wave equations}
\end{center}

Consider at first, the linear problem $\left( 1.7\right) -\left( 1.9\right) $%
. From Theorem 4.2 we obtain the following result

\textbf{Theorem 8.1. } Suppose the the following conditions are satisfied:

(1)\ $a$ is positive, $b$ is a real-valued functions on $\left( 0,1\right) $%
. Moreover$,$ $a\left( .\right) \in C\left( 0,1\right) $ and%
\[
\exp \left( -\dint\limits_{\frac{1}{2}}^{x}b\left( t\right) a^{-1}\left(
t\right) dt\right) \in L_{1}\left( 0,1\right) ; 
\]

(2) assume the Conditions 4.1 are hold and 
\[
\left\vert \alpha _{k}+\beta _{k}\right\vert
>0,\dsum\limits_{k,j=1}^{m}\alpha _{k}\beta _{j}\neq 0,\text{ }0\leq \alpha
<1,
\]%
\[
\frac{1}{q}+\frac{n}{r}=\frac{n}{2}-\gamma =\frac{1}{\tilde{q}}+\frac{n}{%
\tilde{r}}-2;
\]

Let $n\geq 1$ and   
\[
\varphi \in \mathring{W}^{\gamma ,2}\left( R^{n};W^{2,p}\left( 0,1\right)
\right) ,\text{ }\psi \in \mathring{W}^{\gamma -1,2}\left(
R^{n};W^{2,p}\left( 0,1\right) \right) ,\text{ }
\]%
\[
F\in L^{\tilde{q}^{\prime }}\left( \left[ 0,T\right] ;L^{\tilde{r}^{\prime
}}\left( R^{n};L^{2}\left( 0;1\right) \right) \right) .
\]
Let $u$ : $\left[ 0,T\right] \times R^{n}\rightarrow L^{2}\left( 0,1\right) $
be a solution to $\left( 1.7\right) -\left( 1.9\right) $. Then%
\[
\left\Vert A^{\alpha }u\right\Vert _{L^{q}\left( \left[ 0,T\right]
;L^{r}R^{n};L^{2}\left( 0,1\right) \right) }+\left\Vert A^{\alpha
}u\right\Vert _{C\left( \left[ 0,T\right] ;L^{2}\left( R^{n};L^{2}\left(
0,1\right) \right) \right) }+
\]%
\[
\left\Vert A^{\alpha }\partial _{t}u\right\Vert _{C\left( \left[ 0,T\right] ;%
\mathring{W}^{2,\gamma -1}R^{n};L^{2}\left( 0,1\right) \right) }\lesssim
\left\Vert A\varphi \right\Vert _{\mathring{W}^{2,\gamma }\left(
R^{n};L^{2}\left( 0,1\right) \right) }+
\]%
\[
\left\Vert A\psi \right\Vert _{\mathring{W}^{2,\gamma -1}\left(
R^{n};L^{2}\left( 0,1\right) \right) }+\left\Vert F\right\Vert _{L^{\tilde{q}%
^{\prime }}\left( \left[ 0,T\right] ;L^{\tilde{r}^{\prime
}}R^{n};L^{2}\left( 0,1\right) \right) }.
\]

\ \textbf{Proof.} Let $H=L^{2}\left( 0,1\right) $ and $A$ is a operator
defined by $\left( 4.1\right) .$ Then the problem $\left( 1.7\right) -\left(
1.9\right) $ can be rewritten as the problem $\left( 1.2\right) $. By virtue
of $\left[ \text{12, 23}\right] $ the operator $A$ generates analytic
semigroup in $L^{2}\left( 0,1\right) $. Hence, by virtue of (1)-(2) all
conditons of Theorem 4.2 are satisfied. Then Theorem 4.2 implies the
assertion.

Consider now, the problem $\left( 1.10\right) -\left( 1.8\right) -\left(
1.9\right) $. In this section, from Theorem 5.1 we obtain the following
result:

\bigskip \textbf{Theorem 8.2. } Suppose the the following conditions are
satisfied:

(1)\ $a$ is positive, $b$ is a real-valued functions on $\left( 0,1\right) $%
. Moreover$,$ $a\left( .\right) \in C\left( 0,1\right) $ and%
\[
\exp \left( -\dint\limits_{\frac{1}{2}}^{x}b\left( t\right) a^{-1}\left(
t\right) dt\right) \in L_{1}\left( 0,1\right) ; 
\]

(2) the function $F:$ $W^{\frac{1}{2},2}\left( 0,1\right) \rightarrow
L^{2}\left( 0,1\right) $ is continuously differentiable and obeys the power
type estimates

\[
F\left( u\right) =O\left( \left\Vert u\right\Vert _{L^{2}\left( 0,1\right)
}^{k}\right) ,\text{ }\left\Vert u\right\Vert _{L^{2}\left( 0,1\right)
}\left\Vert F_{u}\left( u\right) \right\Vert _{L^{2}\left( 0,1\right) }\sim
\left\Vert F\left( u\right) \right\Vert _{L^{2}\left( 0,1\right) }\text{ } 
\]%
for some $k>1;$

(3) $n\geq 4$, $\gamma =\frac{n-3}{2\left( n-1\right) }$, $k=\frac{\left(
n+1\right) ^{2}}{\left( n-1\right) ^{2}+4};$

(4) $\varphi \in \mathring{W}^{2,\gamma }\left( R^{n};W^{2,2}\left(
0,1\right) \right) $, $\psi \in \mathring{W}^{2,\gamma -1}\left(
R^{n};W^{2,2}\left( 0,1\right) \right) ;$

(5) assume the Conditions 4.1 are hold and 
\[
\left\vert \alpha _{k}+\beta _{k}\right\vert
>0,\dsum\limits_{k,j=1}^{m}\alpha _{k}\beta _{j}\neq 0, 
\]%
\[
\frac{1}{q}+\frac{n}{r}=\frac{n}{2}-\gamma =\frac{1}{\tilde{q}}+\frac{n}{%
\tilde{r}}-2. 
\]

Then for $k\geq k_{0}$ there is a $T>0$ depending only on 
\[
\left\Vert \varphi \right\Vert _{\mathring{W}^{2,\gamma }\left(
R^{n};W^{2,2}\left( 0,1\right) \right) }+\left\Vert \psi \right\Vert _{%
\mathring{W}^{2,\gamma -1}\left( R^{n};W^{2,2}\left( 0,1\right) \right) } 
\]%
and a unique weak solution $u$ to $\left( 1.10\right) -\left( 1.8\right)
-\left( 1.9\right) $ with 
\[
u\in L^{q_{0}}\left( \left[ 0,T\right] ;L^{r_{0}}\left( R^{n};W^{2,2}\left(
0,1\right) \right) \right) , 
\]%
where 
\[
q_{0}=\frac{2\left( n+1\right) }{n-3}\text{, }r_{0}=\frac{2\left(
n^{2}-1\right) }{\left( n^{2}-1\right) +4}. 
\]

In addition, the solution satisfies 
\[
u\in C\left( \left[ 0,T\right] ;\mathring{W}^{2,\gamma }\left(
R^{n};W^{2,2}\left( 0,1\right) \right) \right) \cap C^{1}\left( \left[ 0,T%
\right] ;\mathring{W}^{2,\gamma -1}\left( R^{n};W^{2,2}\left( 0,1\right)
\right) \right) 
\]%
and depends continuously on the data.

\ \textbf{Proof.} Let $H=L^{2}\left( 0,1\right) $ and $A$ is a operator
defined by $\left( 1.4\right) .$ Then the problem $\left( 1.10\right)
-\left( 1.8\right) -\left( 1.9\right) $ can be rewritten as the problem $%
\left( 5.1\right) $. By virtue of $\left[ \text{12, 23}\right] $ the
operator $A$ generates analytic semigroup in $L^{2}\left( 0,1\right) $.
Hence, by virtue of (1)-(5), all conditons of Theorem 5.1 are satisfied.
Then Theorem 5.1 implies the assertion.

\bigskip\ \ \ \ \ \ \ \ \ \ \ \ \ \ \ \ \ \ \ \ \ \ \ \ 

\begin{enumerate}
\item H. Amann, Linear and quasi-linear equations,1, Birkhauser, Basel
1995.\ \ 

\item A. Arosio, Linear second order differential equations in Hilbert
space. The Cauchy problem and asymptotic behaviour for large time, Arch.
Rational Mech. Anal., v 86, (2), 1984, 147-180.

\item A. Ashyralyev, N. Aggez, Nonlocal boundary value hyperbolic problems
Involving Integral conditions, Bound.Value Probl., 2014 V. 2014:214.

\item H. Bahouri, P. G%
\'{}%
erard, High frequency approximation of solutions to critical nonlinear wave
equations, Amer. J. Math. 121 (1999), 131--175.

\item J. Bergh and J. Lofstrom, Interpolation spaces: An introduction,
Springer-Verlag, New York, 1976.

\item M. Christ and M. Weinstein, Dispersion of small amplitude solutions of
the generalized Korteweg-de Vries equation. J. Funct. Anal. 100 (1991),
87-109.

\item G. Da Prato and E. Giusti, A characterization on generators of
abstract cosine functions,\ Boll. Del. Unione Mat., (22)1967, 367--362.

\item Denk R., Hieber M., Pr\"{u}ss J., $R$-boundedness, Fourier multipliers
and problems of elliptic and parabolic type, Mem. Amer. Math. Soc. 166
(2003), n.788.

\item E. B. Davies and M. M. Pang, The Cauchy problem and a generalization
of the Hille-Yosida theorem, Proc. London Math. Soc. (55) 1987, 181-208.

\item D. Foschi, Inhomogeneous Strichartz estimates, Jour. Hyperbolic Diff.
Eqs 2, No. 1 (2005) 1--24.

\item H. O. Fattorini, Second order linear differential equations in Banach
spaces, in North Holland Mathematics Studies,\ V. 108, North-Holland,
Amsterdam, 1985.

\item A. Favini, G. R. Goldstein, J. A. Goldstein and S. Romanelli,
Degenerate second order differential operators generating analytic
semigroups in $L_{p}$ and $W^{1,p}$, Math. Nachr. 238 (2002), 78 --102.

\item A. Favini, V. Shakhmurov, Y. Yakubov, Regular boundary value problems
for complete second order elliptic differential-operator equations in UMD
Banach spaces, Semigroup form, v. 79 (1), 2009.

\item J. Ginibre, G. Velo, Generalized Strichartz inequalities for the wave
equation, Jour. Func. Anal., 133 (1995), 50--68.

\item J. A. Goldstein, Semigroup of linear operators and applications,
Oxford, 1985.

\item V. I. Gorbachuk and M. L. Gorbachuk M, Boundary value problems \ for
differential-operator equations, Naukova Dumka, Kiev, 1984.

\item M. Grillakis, Regularity for the wave equation with a critical
nonlinearity, Comm. Pure Appl. Math., 45 (1992), 749--774.

\item S. G. Krein, Linear differential equations in Banach space,
Providence, 1971.

\item C. Kenig, G. Ponce, L. Vega, A bilinear estimate with applications to
the KdV equation, J. Amer. Math. Soc. 9 (1996), 573--603.

\item L. Kapitanski, Weak and yet weaker solutions of semilinear wave
equations, Comm. Part. Diff. Eq., 19 (1994), 1629--1676.

\item S. Klainerman, M. Machedon, On the regularity properties of a model
problem related to wave maps, Duke Math. J. 87 (1997), 553--589.

\item M. Keel and T. Tao, Endpoint Strichartz estimates. Amer. J. Math. 120
(1998), 955-980.

\item V. Keyantuo, M. Warma, The wave equation with Wentzell--Robin boundary
conditions on $L_{p}$-spaces, J. Differential Equations 229 (2006) 680--697.

\item H. Lindblad, C. D. Sogge, Restriction theorems and semilinear
Klein-Gordon equations in (1+3) dimensions, Duke Math. J. 85 (1996), no 1,
227--252.

\item A. Lunardi, Analytic semigroups and optimal regularity in parabolic
problems, Birkhauser, 2003.

\item M. Meyries, M. Veraar, Pointwise multiplication on vector-valued
function spaces with power weights, J. Fourier Anal. Appl. 21 (2015)(1),
95--136.

\item N. Masmoudi, K. Nakanishi, From nonlinear Klein-Gordon equation to a
system of coupled nonlinear Schrdinger equations, Math. Ann. 324 (2002), (2)
359--389.

\item S. Piskarev and S.-Y. Shaw, Multiplicative perturbations of semigroups
and applications to step responses and cumulative outputs, J. Funct. Anal.
128 (1995), 315-340.

\item A. Pazy, Semigroups of linear operators and applications to partial
differential equations. Springer, Berlin, 1983.

\item C. D. Sogge, Lectures on Nonlinear Wave Equations, International
Press, Cambridge, MA, 1995.

\item C. D. Sogge, Fourier Integrals in Classical Analysis, Cambridge
University Press, 1993.

\item V. B. Shakhmurov, Nonlinear abstract boundary value problems in
vector-valued function spaces and applications, Nonlinear Anal-Theor., v.
67(3) 2006, 745-762.

\item R. Shahmurov, On strong solutions of a Robin problem modeling heat
conduction in materials with corroded boundary, Nonlinear Anal., Real World
Appl., v.13, (1), 2011, 441-451.

\item R. Shahmurov, Solution of the Dirichlet and Neumann problems for a
modified Helmholtz equation in Besov spaces on an annuals, J. Differential
equations, v. 249(3), 2010, 526-550.

\item I. E. Segal, Space-time decay for solutions of wave equations, Adv.
Math. 22 (1976), 304--311.

\item E. M. Stein, Singular Integrals and differentiability properties of
functions, Princeton Univ. Press, Princeton. NJ, 1970.

\item R. S. Strichartz, Restriction of Fourier transform to quadratic
surfaces and decay of solutions of wave equations, Duke Math. J. 44 (1977),
705--774.

\item H. Triebel, Interpolation theory, Function spaces, Differential
operators, North-Holland, Amsterdam, 1978.

\item T. Tao, \ Local and global analysis of nonlinear dispersive and wave
equations, CBMS regional conference series in mathematics, 2006.

\item T. Tao, Low regularity semi-linear wave equations. Comm. Partial
Differential Equations, 24 (1999), no. 3-4, 599--629.

\item C. Travis and G. F. Webb, Second order differential equations in
Banach spaces, Nonlinear Equations in Abstract Spaces, (ed. by V.
Lakshmikantham), Academic Press, 1978, 331-361.

\item G. B. Whitham, Linear and Nonlinear Waves, Wiley--Interscience, New
York, 1975.

\item S. Yakubov and Ya. Yakubov, Differential-operator Equations. Ordinary
and Partial \ Differential Equations , Chapman and Hall /CRC, Boca Raton,
2000.
\end{enumerate}

\bigskip

\end{document}